\documentclass{article}

\usepackage{amsmath,amsfonts,amsbsy,amssymb,fullpage,makeidx,array,graphicx}
\usepackage{microtype}
\usepackage{amsmath}
\usepackage{amssymb}
\usepackage{mathtools}
\usepackage{rotating}
\usepackage{epstopdf}
\usepackage{caption}
\usepackage{lipsum}
\usepackage{afterpage}
\usepackage{subfigure}
\usepackage{pstricks}
\usepackage{bm}
\usepackage[normalem]{ulem}
\usepackage[nomarkers,nolists]{endfloat}
\title{Non-smooth Bifurcations of Mean Field Systems of Two-Dimensional Integrate and Fire Neurons}
\author{Wilten Nicola and Sue Ann Campbell}
\date{\today} 

\begin{document}
\maketitle
\tableofcontents
\begin{abstract}
Mean-field systems have been recently derived that adequately predict the behaviors of large networks of coupled integrate-and-fire neurons \cite{us2}.  The mean-field system for a network of 
neurons with spike frequency adaptation is typically a pair of differential equations for the mean adaptation and mean synaptic gating variable of the network.  These differential equations are non-smooth, and in particular are piecewise smooth continuous (PWSC).  Here, we analyze the smooth and non-smooth bifurcation structure of these equations and show that the system is organized around a pair of co-dimension two bifurcations that involve, respectively, the collision between a Hopf equilibrium point and a switching manifold, and a saddle-node equilibrium point and a switching manifold.  These two co-dimension 2 bifurcations can coalesce into a co-dimension 3 non-smooth bifurcation.  As the mean-field system we study is a non-generic piecewise smooth continuous system, we discuss possible regularizations of this system and how the bifurcations which occur are related to non-smooth bifurcations displayed by generic PWSC systems.  
\end{abstract}
\section{Introduction}

Recently, a class of two-dimensional integrate and fire type models have been developed which can be fit to properties of real neurons.  This class of models includes the adapting, leaky integrate and fire neuron (LIF) \cite{VreeswijkandHansel} , the Izhikevich model \cite{Izhikevich}, the quartic integrate and fire model \cite{Touboul2008}, and the adaptive exponential integrate and fire model (AdEx) \cite{AdEx,Ad2}.  The models in this class are far simpler to fit and simulate than traditional conductance based models. Nevertheless, these models still replicate the more complex behaviors of real neurons \cite{Izhikevich}. These models have been fit to several different neuron types so that the behavior of large networks of these neuron models may be studied through numerical exploration of the parameter space. For example, this approach has been used to determine the role of various parameters in the generation of adaptation induced bursting in networks of CA3 pyramidal 
neurons \cite{Katie,us2}
While the numerical simulation of integrate and fire networks is far simpler and faster than that of conductance based models, numerical exploration of the parameter space is still a time-consuming process.  Furthermore, one cannot easily perform direct bifurcation analysis on large networks.  

Fortunately, a system of mean field equations has been derived for these large networks of two-dimensional integrate and fire neurons \cite{us2}.   This derivation assumes that the networks are all-to-all coupled and the neuronal parameters are homogeneous within the network.  The resulting mean field system is a set of non-smooth differential equations governing the first moments of the adaptation variable and the synaptic coupling variable. 
The mean field system of equations is analytically derived from the original network, without any further fitting.   Thus, one can conduct bifurcation analysis (either analytically or numerically) on the mean field system with confidence that the results are representative of the behavior of the original network of model neurons,  and possibly the original network of actual neurons.  Such a level of correspondence between the parameters for the individual neurons, and the resulting behavior of the full network is currently not possible with more sophisticated types of neuron models. 
 
However, analysis of the derived mean field system has an added level of difficulty as the system of equations is non-smooth.  Both classical bifurcation theory and the newer field of non-smooth bifurcation theory must be used to adequately understand the behavior of the mean field system, and thus the full network.   

Here we explore, both analytically and numerically, many of the non-smooth bifurcations and phenomena that occur in the mean field system of equations from \cite{us2}.  
The primary mean field system we use is that of the Izhikevich model, with the neuronal models fit to hippocampal area CA3 pyramidal neuron data \cite{us}.  We modify the parameters slightly as the neuronal model used in \cite{us} was an alteration of the default Izhikevich model to better fit the action potential half-width observed in the data.     We use this model primarily for two reasons: it is the most analytically tractable and the parameters have been fit to neuronal data.   However,  as we will see, many of the non-smooth bifurcations are present in the other models in the general class of two-dimensional adapting integrate and fire neurons.   Whenever possible we present our results in terms of this general class.

\section{The Mean Field System}
\subsection{The Full Network} 
We consider two-dimensional integrate and fire models of the form 
\begin{eqnarray}
\dot{v} &=& F(v) - w + I \label{de1} \\ 
\dot{w} &=& a(bv - w) \label{de2},
\end{eqnarray}
where $v$ represents the nondimensionalized membrane potential, and $w$ serves as an adaptation variable.  Time has also been non-dimensionalized. 
The dynamical equations (\ref{de1})-(\ref{de2}) are supplemented by the 
following discontinuities 
\begin{equation}
v(t_{spike}^-) = v_{peak} \rightarrow 
\begin{array}{rcl}
v(t_{spike}^+) &=& v_{reset} \\
w(t_{spike}^+) &=& w(t_{spike}^-) + \hat{w}.
\end{array}
\end{equation}
This particular notation was formally introduced in \cite{Touboul2008}, along with a full bifurcation analysis of this general family of adapting integrate and fire neurons.  Members of this family include the Izhikevich model, the adaptive exponential (AdEx) integrate and fire model, and the quartic integrate and fire model \cite{Touboul2008}.

These neurons can be coupled together via a synaptic gating variable, $s(t)$.  The gating variable typically takes the form 
\begin{equation}
s_{ij}(t) = \sum_{t_{j,k}<t}E(t-t_{j,k}).
\end{equation}
where $t_{j,k}$ is the time that the $j$th neuron fires its $k$th spike, $j$ is the index of the presynaptic neuron, and $i$ is the index of the postsynaptic neuron.  The function $E(t)$ varies depending on which synaptic pulse function is used.  Examples include the exponential synapse, the double exponential synapse, and the alpha synapse \cite{us2,bardbook}.  For simplicity we restrict our attention to the exponential synapse, however the analysis can be extended to the other synaptic types without much difficulty.  For the exponential synapse, $E(t)$ is given by :
$$E(t) = \hat{s}\exp\left(\frac{-t}{\tau_s}\right)$$ 

Given the form for $E(t)$, one can
derive a differential equation for $s_i(t) = \sum_{j=1}^{N} s_{ij}(t)$, the total synaptic input to the $i^{th}$ neuron \cite{bardbook,us2}.  For example, for the exponential synapse the differential equation for $s_i(t)$  is 
\begin{eqnarray}
\frac{ds_i(t)}{dt} = -\frac{s_i}{\tau_s} + \frac{\hat{s}}{N}  \sum_{j=1}^N\sum_{t_{j,k}<t} \delta(t-t_{j,k}).
\end{eqnarray}
For all-to-all coupling, the function $s_i(t)$ becomes identical for all the neurons, and can be replaced by a single variable $s(t)$, the global synaptic coupling function.   In this case, the equations for the entire network are: 
\begin{eqnarray}
\dot{v}_i &=& F(v_i) - w + I + {g}s(t)(e_r - v_i) \label{v_fullnet}\\
\dot{w}_i &=& a(bv_i - w_i)\\
\dot{s}&=& -\frac{s}{\tau_s} + \frac{\hat{s}}{N}  \sum_{j=1}^N\sum_{t_{j,k}<t} \delta(t-t_{j,k})\\
v_i(t_{spike}^-) &=& v_{peak} \rightarrow 
\begin{array}{rcl}
v_i(t_{spike}^+) &=& v_{reset} \\
w_i(t_{spike}^+) &=& w_i(t_{spike}^-) + \hat{w}, \label{wjump_fullnet}
\end{array}
\end{eqnarray}
The specific forms of $F(v)$ we consider are:
\begin{eqnarray*}
F(v) &=& -\frac{v}{\tau_m} \quad \text{(Leaky Integrate and Fire)} \\
F(v) &=& v(v-\alpha) \quad \text{(Izhikevich Model)} \\
F(v) &=& e^v - v \quad \text{(Adaptive Exponential Model)}\\
F(v) &=& v^4 -2av \quad \text{(Quartic Model)} 
\end{eqnarray*}
These forms can be arrived at through a suitable non-dimensionalization of the original equations for these models \cite{Touboul2008}.  Note that the non-dimensionalization for the Izhikevich model differs from the one used by \cite{Touboul2008} and is from \cite{us2}.  

\begin{figure}
\centering
          \subfigure[$I=0.4260$, $g=1.2308$]{\includegraphics[scale=0.5]{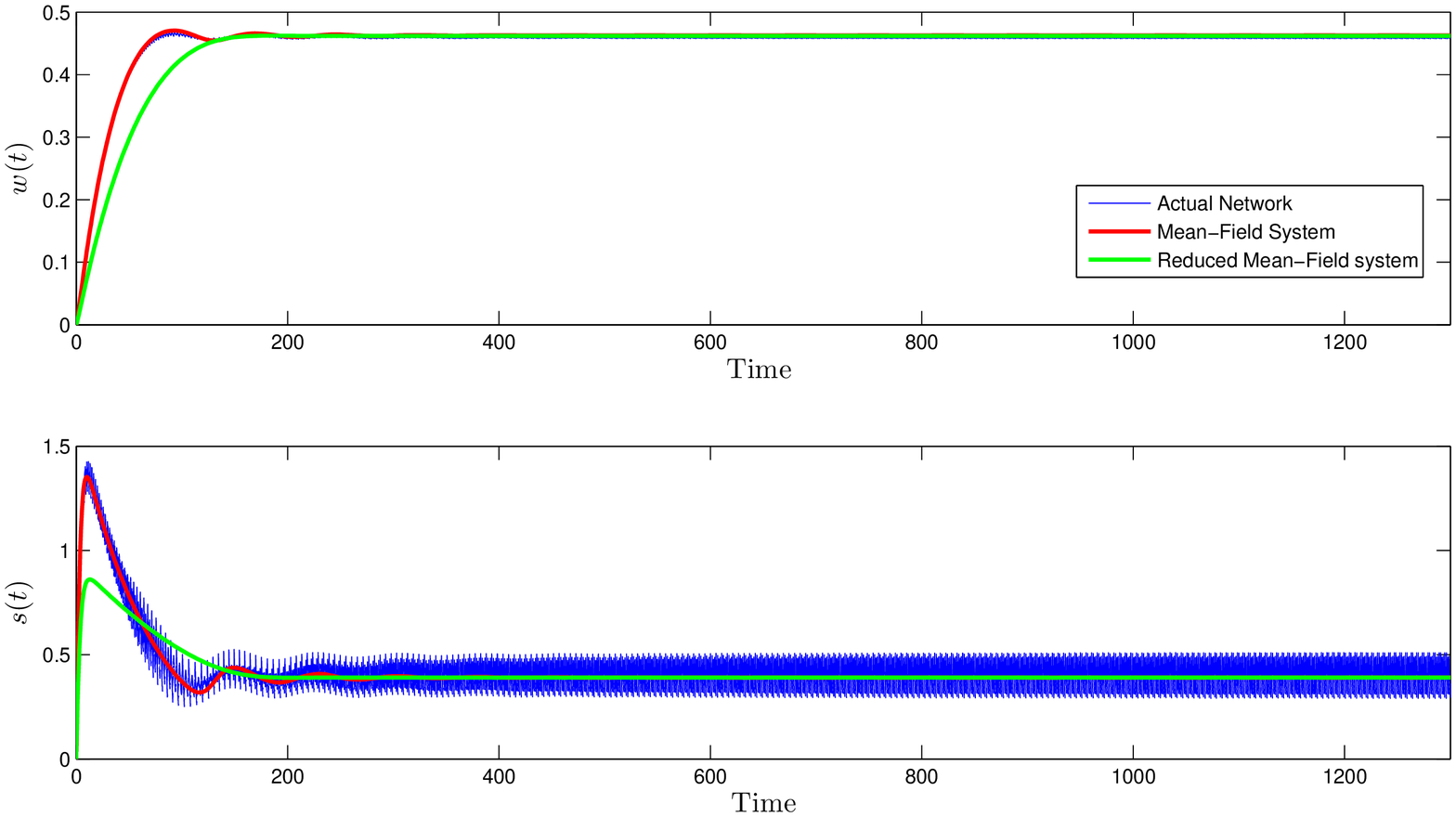}
\label{fig1a}}
         \\
           \subfigure[$I=0.1893$, $g=1.2308$]{\includegraphics[scale=0.5]{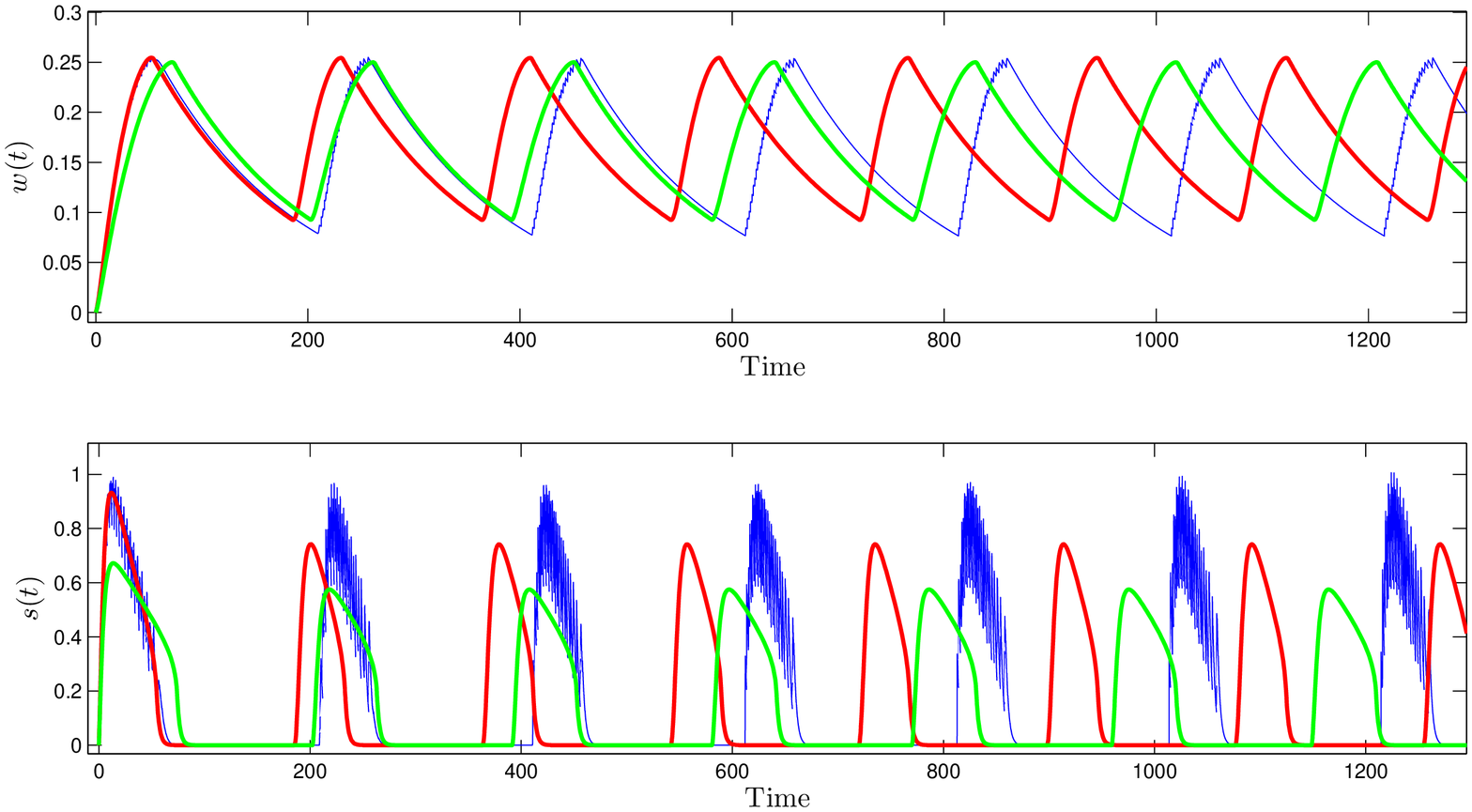}
\label{fig1b}}
        \caption{Actual Network versus mean field Equations.  Simulations of a network of 1000 all-to-all coupled Izhikevich neurons (blue) plotted together with simulations of the mean field equations given by equations (\ref{MFa})-(\ref{MFc}) (red) and of the mean field system resulting from using the simplified firing rate given by equation (\ref{MFs}) (green).  Both mean field models are good approximations to the full system in either tonic firing (\ref{fig1a}) or bursting regimes (\ref{fig1b}).} \label{fig1}
\end{figure}

These networks often display bursting, a oscillatory behaviour where the individual neurons alternate between firing and quiescence \cite{us,us2}.   The other common behaviour is tonic firing, where the neurons all fire at a constant rate. The transition between these two behaviours is a bifurcation of the full network.  An example of this transition for a network of 1000 Izhikevich neurons with all-to-all coupling is shown in Figure \ref{fig1}.  The parameters for this network can be found in Table 1.  In Figure \ref{fig1a}, the neurons in the network fire spikes, and the mean-adaptation variable, $w $, and the synaptic coupling variable, $s$ both converge to a stable steady state.  In Figure \ref{fig1b}, the neurons fire synchronized bursts, and the pair of variables $(w ,s)$ converge to a steady state limit cycle, representing the oscillation between firing and quiescence that the individual neurons undergo.  This occurs as the current $I$ is decreased from $I=0.4260$ in Figure \ref{fig1a} to $I=0.1893$ in Figure \ref{fig1b}.  

As seen in Figure \ref{fig1}, the network variables $s$ and $w $ can capture a great deal of information about the behavior of the entire network.  Further, their steady state behaviour undergoes a qualitative change as the parameter $I$ is decreased.  Thus, it would be advantageous to have a closed set of differential equations for these variables, as any qualitative change in the behavior of the full network should manifest itself as a bifurcation of the dynamical system for $(s,w )$.  In the following subsection, we derive such a system for these network variables.  

\subsection{The Mean Field System} 

To derive the mean field system one begins by defining the population 
density function, $\rho(v,w,t)$, which is  a probability density function
for the location of of the variables $v,w$ in the phase space. That is,
the probability of finding  a neuron in the region $\Omega$ of phase space  
is given by integrating $\rho$ over $\Omega$.
Starting from the full network model, one can derive (in the limit that 
$N\rightarrow\infty$) a partial differential equation that governs the 
evolution of the probability density function $\rho(v,w,t)$ and predicts
the large network dynamics of the original model \cite{NT,OKS}.
This partial differential equation, called the population density equation, 
takes the form 
\begin{eqnarray*}
\frac{\partial \rho(v,w,t)}{\partial t}  &=& -\nabla \cdot \bm{J}(v,w,t) \\
\bm{J}(v,w,t) &=& \begin{pmatrix} J_V(v,w,t) \\ J_W(v,w,t) \end{pmatrix} = \rho(v,w,t)\begin{pmatrix}F(v) - w + I + gs(e_r-v) \\ a(bv-w) \end{pmatrix}
\end{eqnarray*}
where the term $\bm{J}(v,w,t)$ is referred to as the flux.  The discontinuities and discrete jumps in the model neurons impose a boundary condition on the flux:
\begin{eqnarray*}
J_V(v_{peak},w,t) = J_V(v_{reset},w+\hat{w},t).
\end{eqnarray*}
This PDE is coupled to an ODE for $s$, given by 
\begin{equation}
s' = -\frac{s}{\tau_s} + \hat{s}\int_W J_V(v_{peak},w,t)\,dw 
\end{equation}

In order to reduce this system to a small, closed set of ordinary differential equations, one has to apply a series of approximations.  The derivation is somewhat lengthy thus we refer the reader to \cite{us2} for the exact details.   After the approximations are made, the resulting mean field system is given by: 

\begin{eqnarray}
s' &=& -\frac{s}{\tau_s} + \hat{s}\langle R_i(t)\rangle \label{MFa} \\
w' &=& -\frac{w}{\tau_w} + \hat{w}\langle R_i(t)\rangle \\
\langle R_i(t) \rangle &=&\begin{cases} \left[ \int_{V}\frac{dv}{F(v)-w  + I + g(e_r - v)s}  \right]^{-1} & \quad \text{if}\quad  H(w  ,s)>0\\ 0 & \quad \text{if} \quad H(w ,s)\leq 0  \end{cases} \label{MFc}
\end{eqnarray}
Here $s$ and $w $ correspond to the mean network adaptation and global synaptic coupling variable.  Note that we have omitted the $\langle \rangle$ brackets denoting the average value of $w$ present in \cite{us2} for simplicity and clarity. The function $\langle R_i(t)\rangle$ is the 
instantaneous network averaged firing rate, as a function of $s$ and $w $. The function $H$ defines when the integral in \eqref{MFc} makes 
sense. It defines {\em switching manifold} of the nonsmooth system
\eqref{MFa}--\eqref{MFc}.

One can derive an expression for the switching manifold equation by determining
when the denominator in \eqref{MFc} 
\[
F(v) - w  + I + gs(e_r - v) \stackrel{def}{=} G(v,s,w )
\]
first becomes zero somewhere in the $s,w$ phase space.  To do this we find
the minimum of $G(v,s,w )$, for $v\in[v_{reset},v_{peak}]$, 
regarding $s$ and $w $ as fixed parameters.   
For the general class of models studied in \cite{Touboul2008}, the function $F(v)$ is 
assumed to be strictly convex, that is $F''(v)>0$, and at least three times 
continuously differentiable.  It follows that $G(v,s,w )$ is also strictly convex as a function of $v$, and hence its minimum on 
$[v_{reset},v_{peak}]$ occurs at a critical point.  The critical points 
(as a function of $v$) are given by solving 
\begin{eqnarray}
\frac{\partial G}{\partial v} &=& F'(v) -gs =  0\Rightarrow  \nonumber \\
F'(v^*(s)) &=& gs \label{exp0}
\end{eqnarray}
Thus $v^*(s)$ is the value of $v$ at which $G$ has a minimum.  

The leaky integrate-and-fire neuron is not strictly convex. In fact
$F'(v)<0$ for this model, thus the minimum of $G$ occurs at the right endpoint 
of the interval.

In either case, the minimum value defines the function $H$
\begin{equation}
H(s,w )=G(v^*(s),s,w ) = F(v^*(s)) -w  + gs(e_r - v^*(s)) + I, 
\label{Hdef}
\end{equation}
and the switching manifold equation 
$$0=H(s,w )= I - w   + F(v^*(s)) + gs(e_r - v^*(s)) = I-I^*(s,w ).$$ 
This latter expression is useful as we can think of $I^*$ as an $s-$ and $w-$
dependent rheobase current. Anywhere in the phase space where $I-I^*(s,w)>0$ the
network is firing with mean firing rate given by 
\begin{equation}
\langle R_i(t)\rangle = \left[ \int_{V}\frac{dv}{F(v)-w  + I + g(e_r - v)s}  \right]^{-1}.  \label{fr}
\end{equation}
Anywhere that $I-I^*(s,w)\le 0$ the network is quiescent and 
the mean firing rate is 0.

There are a couple of important facts to note before we proceed further.  
First of all, $I^*(0,0)=I_{rh}$, the rheobase current for the uncoupled, 
nonadapting neuron, which is governed by the equation 
$$\dot{v} = F(v) + I .$$
Based on the assumptions made on $F(v)$ in \cite{Touboul2008}, then this model neuron has a type-I firing profile.  
Additionally, given that $I^*(0,0) = -F(v^*(0))$, we also have the 
following equation
$$F(v^*(0)) = -I_{rh}.$$ 
These two facts will prove important for our later analysis.

To conclude we display some expressions for specific models.
The rheobase currents are given by 
\begin{eqnarray*}
I^*(s,w ) &=& w  - gse_r + v_{peak}\left(\frac{1}{\tau_m} + gs\right)\quad \text{(Leaky Integrate and Fir)}\\
I^*(s,w ) &=&  w  - gse_r + \frac{(\alpha+gs)^2}{4} \quad \text{(Izhikevich)}\\
I^*(s,w ) &=&  w  - gse_r + (1+gs)(\log(1+gs)-1) \quad \text{(AdEx)}\\
I^*(s,w ) &=&  w  - gse_r +3\left(\frac{gs+2a}{4}\right)^{\frac{4}{3}} \quad \text{(Quartic)}\\
\end{eqnarray*}
with corresponding minimum values of $G$:
\begin{eqnarray*}
v^*(s) &=& v_{peak}\\
v^*(s) &=& \frac{\alpha+gs}{2}\\
v^*(s) &=& \log(1+gs)\\
v^*(s) &=& \left(\frac{gs+2a}{4}\right)^{1/3}\\
\end{eqnarray*}

For both the LIF model and Izhikevich model, the mean firing rate 
can be evaluated analytically: 
\begin{eqnarray*}
 \langle R_i(t)\rangle &=& \frac{\frac{1}{\tau_m} + gs}{\log \left(\frac{-v_r\left(\frac{1}{\tau_m }+ gs\right)+I+gse_r - w}{-v_p\left(\frac{1}{\tau_m }+ gs\right)+I+gse_r - w} \right)} \quad \text{(LIF firing rate)} \\
\langle R_i(t)\rangle &=& \frac{\sqrt{I-I^*(s,w )}}{\arctan\left( \frac{v_{peak} - \frac{1}{2}(\alpha+gs)}{\sqrt{I-I^*(s,w )}}\right)-\arctan\left( \frac{v_{reset} - \frac{1}{2}(\alpha+gs)}{\sqrt{I-I^*(s,w )}}\right)} \quad \text{(Izhikevich Firing Rate)} 
\end{eqnarray*}
For the other models, the firing rate must be evaluated numerically. This can be done by integrating equation \eqref{fr} over $[v_{reset},v_{peak}]$ treating $ w $, and $s$ as fixed parameters at each time step. 
This approach can be used to numerically analyze the bifurcation types of these equations using numerical bifurcation software, such as MATCONT, among others \cite{MATCONT}.  However, the numerical integration method should be of high enough order accuracy for the numerical continuation results to be trusted.   

Given the mean field system described above, one should consider whether numerical bifurcation results or analytical approaches should be taken.  Numerical bifurcation analysis can yield results which are accurate throughout the phase space, but require choosing a particular model and determining which parameters to fix and which to vary.  Analytical methods can yield model independent results and give insight into the role of various parameters in system behaviour, however, they are often restricted to particular regions of the parameter space and/or phase space, as we shall see.  

\section{Analytical Results} 

In order to proceed analytically, we need to sacrifice some of the complexity of the original mean field system.  In particular, as the usual formulas of $\langle R_i(t)\rangle$ are difficult to deal with analytically, we approximate the firing rate for all the two dimensional integrate and fire models as follows:
\begin{eqnarray*}
\langle R_i(t) \rangle \approx f(I-I^*(s, w))
\end{eqnarray*}
where the form of $f$ varies from model to model.   

There are two approaches to justify this particular approximation, one local and one global.  The local method relies on the fact that all the neuron models fire with a type I profile, and as such, they have a characteristic firing rate near the switching manifold that is proportional to 
\begin{equation} \langle R_i(t) \rangle \propto \sqrt{I-I^*(s, w )} \label{fr1} \end{equation}  
This basic result follows from the normal form for a saddle node bifurcation.  This reduction is valid when $I-I^*(s, w )$ is small.  This occurs when $ w$ and $gs$ are small and $I$ is close to $I_{rh}$.    We note that this approximation and very similar differential equations appear in \cite{bardbook,bardpaper}.  For example, equations (3.6)-(3.7) in \cite{bardpaper} are similar to ours however the interpretation for those particular equations was for the firing rate of an E/I coupled pair of neurons.  Additionally, the non-smooth nature of those equations is not explored to the best our knowledge in any source.  

\begin{figure}
\centering
\includegraphics[scale=0.8]{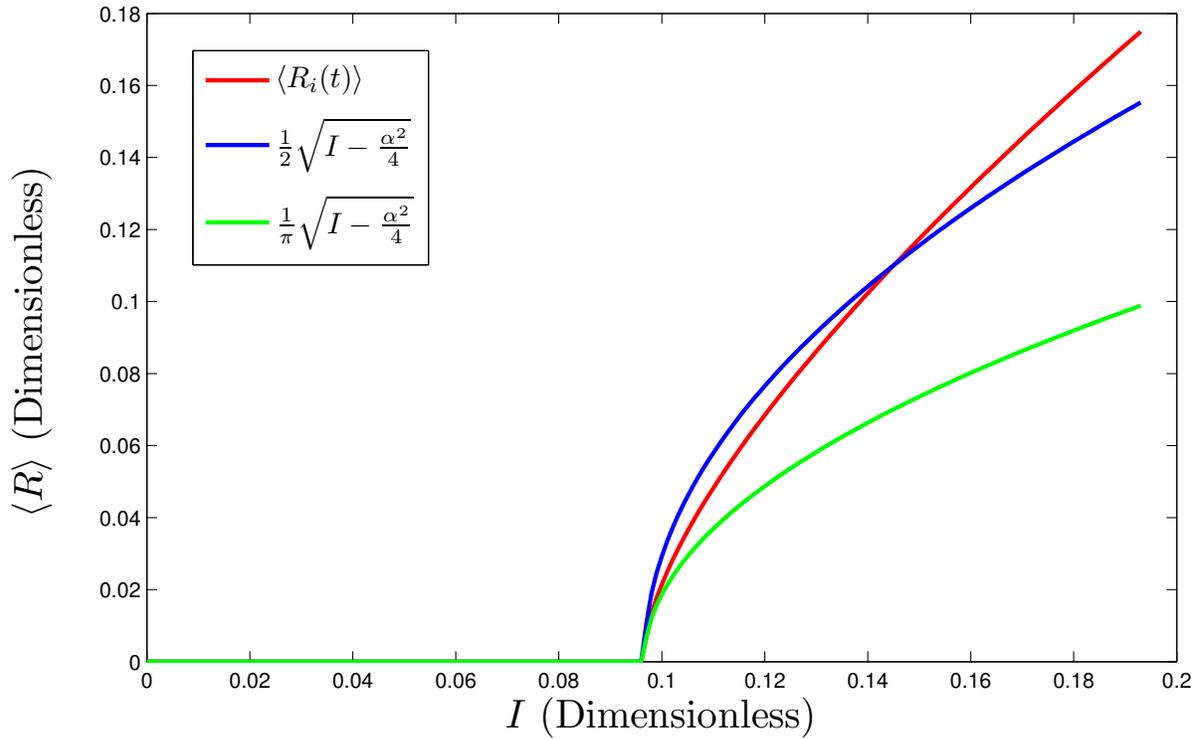}
\caption{The approximations for $\langle R_i(t)\rangle$ (in red) using $k \sqrt{I-I^*(s,w)}$ for $k$ predicted from topological normal form theory ($k = 1/\pi$, green), which is only locally valid near the transition to firing from quiescence, or for a more global fit ($k=1/2$, blue).  The specific $k$ used does not matter as it can be merely absorbed into the $\hat{s}$ and $\hat{w}$ parameters and the final bifurcation analysis is the same.}\label{fig2}
\end{figure}

For a more global approximation, one can fit a function $f(I-I^*(s, w ))$ to $\langle R_i(t)\rangle$.   For example, we have found that for the Izhikevich network, the approximation 
\begin{equation}
\langle R_i(t) \rangle = \left [ \int_{v_{reset}}^{v_{peak}} \frac{dv}{v(v-\alpha)+ I -  w + gs(e_r-v)} \right]^{-1} \approx \frac{1}{2}\sqrt{I - I^*(s, w)} \label{MFs}
\end{equation}
is much better globally than the local approximation for type I firing neurons, which is given by: 
$$ \langle R_i(t) \rangle  \approx \frac{1}{\pi} \sqrt{I-I^*(s,w)}$$ 
as shown in figure \ref{fig2}.  
As equations \eqref{fr1} and \eqref{MFs} only differ by a multiplicative constant, we take the firing rate to be in the form $\sqrt{I-I^*(s, w)}$, which is valid for the Izhikevich network globally, and for the other networks locally in the phase and parameter space.   Note that whichever constant is used, it can be merely absorbed into the $\hat{s}$ and $\hat{w}$ terms in the mean field dynamical equations.  

This particular approximation turns out to yield a system that is tractable to analysis of both the smooth and non-smooth bifurcations displayed by this system, and shows a considerable degree of accuracy with both the actual network, and the original mean-field system, as shown in figure \ref{fig1}.   We reiterate that the non-smooth nature of these equations is often not considered.  However, we will show that recent developments in non-smooth bifurcation theory allow us to analyze these equations more fully.  

With the above simplification, the approximate mean field system that we analyze is given by 
\begin{eqnarray}
\dot{s} &=& f(s,w)= -\frac{s}{\tau_s} + \hat{s}\,\langle R_i(t) \rangle  \label{s1}\\
\dot{w} &=& g(s,w)=-\frac{w}{\tau_w} + {\hat{w}}\,\langle R_i(t) \rangle \label{w1}\\
\langle R_i(t) \rangle &=& \begin{cases} \sqrt{I-I^*(s,w )} & \quad  I-I^*(s,w )> 0\\  0 & \quad  I-I^*(s,w)\leq 0\end{cases} \label{r1}
\end{eqnarray}
where the switching manifold varies depending on which neuron model is used.  The terms $\hat{s}$ and $\hat{w}$ are rescaled to absorb any constant term in (\ref{fr1}) or (\ref{MFs}).   We will refer to equations (\ref{s1})-(\ref{r1}) from here on as the reduced mean-field system.  

We first classify the type of non-smoothness that this system exhibits.  The system is smooth and has derivatives of all orders everywhere except on the switching manifold, i.e., when $I-I^*(s,w )=0$.  On the switching manifold, the system is continuous but not differentiable. Thus, this is a piecewise-smooth continuous (PWSC) system. Equivalently, it has a uniform order of discontinuity of 1 \cite{nonsmooth}.

In contrast, the order of discontinuity for the full (non-simplified) Izhikevich mean field system is non-uniform.  To see this, rewrite $I-I^*(w,s)=h(I,s)-w$  
and consider the limit
\[
\lim_{w\rightarrow h(I,s)}\langle R_i(t)\rangle = 
\lim_{w\rightarrow h(I,s)}
\frac{\sqrt{h(I,s)-w}}{\arctan\left( \frac{v_{peak} - \frac{1}{2}(\alpha+gs)}{\sqrt{h(I,s)-w}}\right)-\arctan\left( \frac{v_{reset} - \frac{1}{2}(\alpha+gs)}{\sqrt{h(I,s)-w}}\right)} 
\]
Straightforward calculations show that this limit is $0$, and hence the
order of discontinuity of the model is $1$, only if
\begin{eqnarray*}
 2v_{reset}-\alpha< gs< 2v_{peak}-\alpha.
\end{eqnarray*}
Outside of this region the order of discontinuity is not $1$. In general
$2v_{reset}-\alpha<0$ for physical reasons, so the order of discontinuity
of the full Izhikevich model is $1$ if $gs$ is sufficiently small.
Thus the approximation simplifies the analysis as it has a uniform order of
discontinuity, however, it will likely be a reasonable approximation of the
Izhikevich model only when this latter model has order of discontinuity $1$,
i.e., if $gs$ is sufficiently small.

To summarize, for all the neuron models we consider the firing rate to be approximated by $\sqrt{I-I^*(s,w)}$ based on the type I profile for all the models, where $I^*(s,w)$ varies depending on which model is used.  These approximations are only local in nature, being valid only when $I\approx I^*(s,w)$, aside from the Izhikevich model, as one can fit a more global approximation to the network averaged firing rate, as stated earlier.  

\subsection{Existence and Linear Stability of Equilibria}\label{stabsec}  
The equilibria of the mean field equations \eqref{s1}-\eqref{w1} depend 
on the sign of $I-I^*(s,w )\le 0$. 

If  $I-I^*(s,w)\le 0$ then the only equilibrium point is 
the trivial solution, $e_0=(0,0)$, which is a stable node.  This equilibrium 
corresponds to all 
the neurons being quiescent,  $\langle R_i(t) \rangle = 0$, thus we will 
refer to it as the non-firing solution. It will only exist when
the origin of the phase space lies in the region where $I-I^*(0,0)\le 0$,  
which corresponds to $I\le I_{rh}$.   Alternatively, in the language of non-smooth dynamical systems theory, $e_0$ is virtual if $I > I_{rh}$ and real if $I\leq I_{rh}$ \cite{nonsmooth}.  

If $I-I^*(s,w )> 0$, nontrivial equilibria $(s,w)$ may
exist. If they do then $s,w$ must satisfy
\begin{eqnarray}
s &=& \lambda_s \sqrt{I-I^*(s,w)}\label{work1} \\
w &=& \lambda_w \sqrt{I-I^*(s,w)}\label{work2}
\end{eqnarray}
where $\lambda_s = \tau_s \hat{s}$, and a similar equation holds for $\lambda_w$.  Equations (\ref{work1}) and (\ref{work2}) yield the following relationship 
\begin{equation}
w = \frac{\lambda_w}{\lambda_s} = \eta s. \label{work4}
\end{equation}
Thus the equilibria are given by $(s,\eta s)$ where $s$ satisfies the
nonlinear equation
\begin{equation}
s = \lambda_s\sqrt{I - I^*(s,\eta s)}. \label{seqm}
\end{equation}
Note that equation \eqref{work1} implies that $s = \lambda_s \langle R_i(t) \rangle \geq 0$.  Thus for an equilibrium to be a valid, it must satisfy $s\geq 0$ (which implies $w\ge 0$). 

The equilibrium condition \eqref{seqm} for the quartic and AdEx models yield nonlinear equations without analytic closed form solutions.   However, one can apply a power series (assuming that $gs$ is small) to come up with an approximation to the steady solutions.   Note that $v^*$ is actually always a function of $gs$, as opposed to just $s$, as it is given by solving the algebraic 
equation \eqref{exp0}.
Thus, we can write down the following expansions for $v^*(s)$ and  $F(v^*(s))$ 
\begin{eqnarray} 
v^*(s) &=& v^*(0) + {v^*}'(0)gs + O((gs)^2) \label{exp1}\\
F(v^*(s)) &=& F(v^*(0)) + {v^*}'(0)\frac{(gs)^2}{2} + O((gs)^3) \nonumber \\ 
               &=&  -I_{rh} + {v^*}'(0) \frac{(gs)^2}{2} \label{exp2}
\end{eqnarray}
where (\ref{exp2}) can be derived from using the relationship (\ref{exp0}).
Since these expansions are valid if $gs$ is sufficiently small, we will refer to them as
the weak coupling expansions.
If we use only the initial terms shown in equations (\ref{exp1})-(\ref{exp2}), then the approximate solution obtained is exact in the case of the Izhikevich model as all higher order terms vanish, and is the simplest analytical approximation to the other models.  Using these terms, we arrive at the following equation for the weak-coupling expansion to the equilibrium points: 
\begin{eqnarray*}
\frac{s^2}{\lambda_s^2} &=& I-I_{rh}  -{v^*}'(0)\frac{(gs)^2}{2} + gs(e_r - v^*(0)) - \eta s \\
0&=& s^2\left(\frac{1}{\lambda_s^2}+{v^*}'(0)\frac{g^2}{2}  \right) + s\left(\eta - g(e_r-v^*(0))\right) + I_{rh}-I \\
0 &=& A_2(g)s^2 + A_1(g)s + A_0.
\end{eqnarray*}
This equation can be solved, yielding two solution branches:
\begin{eqnarray*}
s_{\pm} = -\frac{A_1(g)}{2A_2(g)} \pm \sqrt{\frac{A_1(g)^2}{4A_2(g)^2} - \frac{A_0}{A_2(g)}}.
\end{eqnarray*}
We will denote the corresponding equilibria as 
$e_\pm=(s_\pm,w_\pm)=(s_\pm,\eta s_\pm)$.
Defining the new parameters
\begin{eqnarray}
\tilde{I} &=& -\frac{A_0}{A_2(g)} = \frac{I-I_{rh}}{A_2(g)} \label{tildeIdef}\\ 
\beta &=& -\frac{A_1(g)}{2A_2(g)} =  \frac{g(e_r-v^*(0))-\eta}{2A_2(g)}\label{betadef}\\
&=& \frac{(e_r-v^*(0))}{2A_2(g)}\left(g-\frac{\eta}{e_r-v^*(0)}\right) = M(g)(g-g^*) 
\end{eqnarray}
the $s$ variable of the solution branches may be written
\begin{equation} 
s_{\pm}(\beta,\tilde{I}) = \beta \pm \sqrt{\beta^2 + \tilde{I}}. \label{ss}
\end{equation}
Note that $A_2(g)>0$. Further, since $v^*(0)$ is the minimum of $F(v)$ and the reversal potential for an excitatory synapse is above the resting membrane potential, $v_r$, we have $e_r>v_r > v^*(0)$. 
It follows that $M(g)$ is a strictly positive function.  

In a similar way, one can solve for the equilibrium values of the full Izhikevich model analytically.  Two equilibria $e_{\pm}$ of the same form are obtained 
with $s\pm$ given by: 
\begin{eqnarray*}
s_+(g,I)&=& \frac{-(\eta - g(e_r - \frac{\alpha}{2})) + \sqrt{ (\eta - g(e_r - \frac{\alpha}{2}))^2 + 4(I-\frac{\alpha^2}{4})(\frac{1}{\lambda_s^2} + \frac{g^2}{4})}} { \frac{1}{\lambda_s^2} + \frac{g^2}{4}}\\
s_-(g,I)&=&  \frac{-(\eta - g(e_r - \frac{\alpha}{2})) - \sqrt{ (\eta - g(e_r - \frac{\alpha}{2}))^2 + 4(I-\frac{\alpha^2}{4})(\frac{1}{\lambda_s^2}+ \frac{g^2}{4})}}{\frac{1}{\lambda_s^2} + \frac{g^2}{4}}
\end{eqnarray*}
Introducing the parameters
\begin{eqnarray}
\tilde{I} &=&  \frac{ (I-\frac{\alpha^2}{4})}{\frac{1}{\lambda_s^2}  + \frac{g^2}{4}}   \\ 
\beta &=& -\frac{(\eta - g(e_r - \frac{\alpha}{2}))}{2(\frac{1}{\lambda_s^2} + \frac{g^2}{4})}
\end{eqnarray}
the steady states can again be written in the form \eqref{ss}.
Note, as a check of consistency, that $I_{rh} = \frac{\alpha^2}{4}$, and $v^*(0) = \frac{\alpha}{2}$ for the Izhikevich model.  

\begin{figure}
\centering
          \subfigure[$(\beta,\tilde{I})$ plane ]{\includegraphics[scale=0.8]{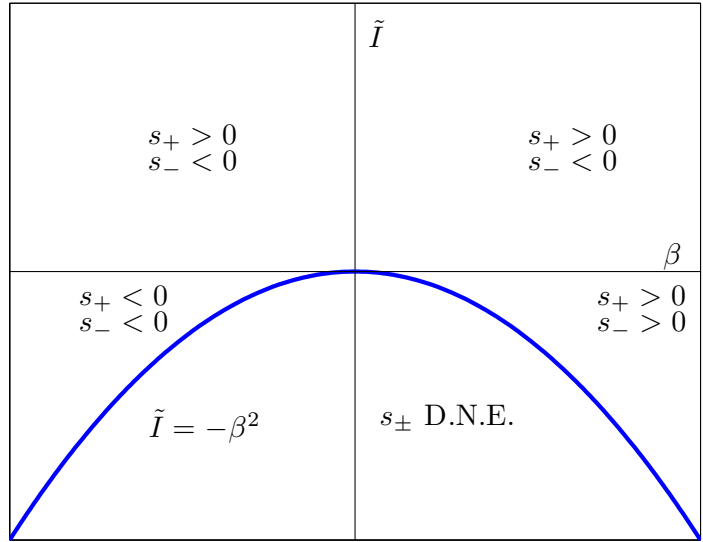}
\label{fig3a}}
\\
           \subfigure[$(g,I)$ plane]{\includegraphics[scale=0.8]{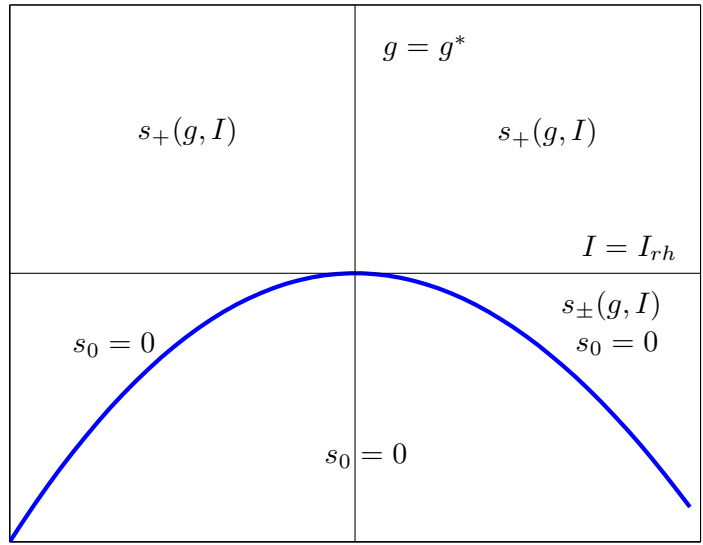}
\label{fig3b}}
        \caption{The existence of equilibria for the mean field system.  
(a) The sign of the $s$ component of the nontrivial equilibria, in the $(\beta,\tilde{I})$ parameter plane.  $s_+>$ is positive in the first two quadrants and in a narrow wedge-shaped region in the fourth quadrant.  $s_-$ is also positive in the wedge-shaped region.  (b) The existence of the trivial and nontrivial
equilibria for the Izhikevich model in the $I,g$ parameters space. The nontrivial equilibrium $e_+(g,I)$ only exists in the regions $I>\alpha^2/4$, and for $I<\alpha^2/4$ in the wedge shaped region of the fourth quadrant indicated.  
The nontrivial equilibrium $e_-(g,I)$ only exists in this wedge shaped region.  
The trivial (non-firing) equilibrium $e_0$ only exists for $I\le \alpha^2/4$.} \label{fig3}
\end{figure}

Based on the form \eqref{ss} and the fact that $A_2(g)>0$, it is 
straightforward to 
show the signs of $s_\pm$ are as shown in Figure \ref{fig3}(a). 
Since we require the equilibrium solutions to be positive, $e_{\pm}$ will 
have different regions of existence depending on the values of $\beta$ 
and $\tilde{I}$. In particular, both equilibrium points exist when 
$I<I_{rh}$ and $g>\frac{\eta}{e_r - v^*(0)}$ in a wedge shaped region given by $\beta^2 + \tilde{I}>0$ .   Only $e_+$ exists when
$I>I_{rh}$.   Neither solution exists in other parts of the parameter space.
The regions of existence of $e_\pm$ and the non-firing solution are
show for the Izhikevich model in Figure \ref{fig3}(b).

Away from the switching manifold, we can analyze the smooth bifurcations 
of the equilibria via linearization.   The non-firing solution does not undergo any smooth bifurcations, as it lies in the region of phase space governed
by the equations
\begin{eqnarray*}
s' &=& -\frac{s}{\tau_s} \\
w'&=& -\frac{w}{\tau_w} 
\end{eqnarray*}
Thus the non-firing solution is asymptotically stable when it exists and
does not lie on the switching manifold, i.e., for $I<I_{rh}$.
To analyze the nontrivial equilibrium points, consider the Jacobian of the 
mean field system in the region of phase space where $I-I^*(s,w)>0$ 
\begin{eqnarray*}
J(s,w) &=& \begin{pmatrix} -\frac{1}{\tau_s} - \frac{\hat{s} \frac{\partial I^*(s,w)}{\partial s}}{2\sqrt{I-I^*(s,w)}} &  -\frac{\hat{s} \frac{\partial I^*(s,w)}{\partial w}}{2\sqrt{I-I^*(s,w)}} \\  -\frac{\hat{w} \frac{\partial I^*(s,w)}{\partial s}}{2\sqrt{I-I^*(s,w)}}& -\frac{1}{\tau_w} - \frac{ \hat{w} \frac{\partial I^*(s,w)}{\partial w}}{2\sqrt{I-I^*(s,w)}}
\end{pmatrix}. 
\end{eqnarray*}
We can simplify the Jacobian by imposing the steady state condition $\sqrt{I-I^*(s,w)} = s/\lambda_s$.  Additionally, one can note that:
\begin{eqnarray*}
\frac{\partial I^*(s,w)}{\partial w} &=& 1 \\
\frac{\partial I^*(s,w)}{\partial s} &=& -g(e_r - v^*(s)). 
\end{eqnarray*}
Applying these formulas yields the following simplified Jacobian 
\begin{eqnarray*}
J(s) &=& \begin{pmatrix} -\frac{1}{\tau_s} + \frac{\lambda_s \hat{s} g(e_r - v^*(s)) }{2s} &  -\frac{\lambda_s \hat{s} }{2s} \\  \frac{\lambda_s \hat{w} g(e_r - v^*(s)) }{2s}& -\frac{1}{\tau_w} - \frac{\lambda_s \hat{w}}{2s}
\end{pmatrix}, 
\end{eqnarray*}
which has trace and determinant given by
\begin{eqnarray*}
tr(J) &=& -\left(\frac{1}{\tau_s} + \frac{1}{\tau_w}\right) +\frac{\lambda_s}{2s}\left(\hat{s} g(e_r - v^*(s))-\hat{w} \right)\\
det(J) &=& \frac{1}{\tau_s\tau_w}\left(1 - \frac{\lambda_s^2}{ 2s}\left(g(e_r-v^*(s)) - \eta\right)   \right).
\end{eqnarray*}
We can now discuss the stability of each equilibrium in its region of existence.  First we apply the weak coupling expansion \eqref{exp1} to the determinant:
\begin{eqnarray*}
\text{det}(J) &=& 
\frac{1}{\tau_s\tau_w}\left(1+\lambda_s^2\frac{g^2}{2}{v^*}'(0) + \frac{\lambda_s^2}{2s}(g(e_r-v^*(0))-\eta)\right)\\
&=& \frac{A_2(g)\lambda_s^2}{\tau_s\tau_w}\left(1+\frac{A_1(g)}{2sA_2(g)}\right).\\
\end{eqnarray*} 
Upon substitution of the equilibrium values of $s$ and using the definition 
\eqref{betadef} of $\beta$, we obtain 
\begin{equation}
\left. \text{det}(J)\right|_{s_\pm} = \frac{A_s(g)\lambda_s^2}{\tau_s\tau_w}\left(1-\frac{\beta}{\beta\pm \sqrt{\beta^2 + \tilde{I}}}\right). \label{detJ}
\end{equation}
Since the sign of $ \frac{A_s(g)\lambda_s^2}{\tau_s\tau_w}$ is strictly positive, and the equilibria are only defined when $\beta \pm \sqrt{\beta^2 + \tilde{I}}\geq 0$, we can immediately conclude that
\begin{eqnarray*}
\left.\text{det}(J)\right|_{s_+}   &\geq& 0 \\
\left.\text{det}(J)\right|_{s_-} &\leq& 0 
\end{eqnarray*}
Now, this implies that the equilibrium $e_-$ is always an unstable saddle. 
The equilibrium $e_+$, however, can be a node or a focus and its stability 
is determined by the trace.  We will discuss this further in section 
\ref{Hopfsec}. Note that these results are only valid when $gs_-$ is small 
for the QIF and AdEx models, but are globally valid for the Izhikevich model.  

We can use the equations for the trace and determinant to formulate necessary conditions for the equilibria to display certain smooth bifurcations.  In particular, $\text{det}(J) = 0$ and $\text{tr}(J) \neq 0$ are necessary conditions for an equilibrium to undergo a saddle-node bifurcation, while $\text{det}(J)>0$ and $\text{tr}(J) = 0$ are necessary conditions for a Hopf bifurcation. Having both $\text{det}(J) = 0$ and $\text{tr}(J)=0$, is a necessary condition for a
Bogdanov-Takens bifurcation. Of course, to determine if these bifurcations 
actually occur requires checking additional genericity conditions.
In the following section, we check these conditions where possible.

\subsection{The Saddle Node Bifurcation Condition} 
As described above, necessary conditions for a saddle-node bifurcation 
are $\det(J)=0$ and $\text{tr}(J) \neq 0$. It is easy to see from \eqref{detJ}
that the first condition is satisfied for both $e_\pm$ when 
$\beta^2 + \tilde{I} =0$. It can be shown that the second condition is
satisfied except at isolated points in the $(g,I)$ parameter space as 
determined in section~\ref{Hopfsec}. In the following we will assume that
we exclude these points.

To pursue this further, we study the existence of the equilibria. From the 
previous subsection, we know that 
$e_\pm$ both exist if $\beta^2 + \tilde{I}>0$ and neither exists if 
$\beta^2 + \tilde{I}<0$.  When $\beta^2 + \tilde{I} = 0$, 
the two equilibria collapse into a single equilibrium, with $s = {\beta}$. 
We thus conclude that $ \tilde{I} =-\beta^2$ corresponds to a two-parameter 
curve of saddle-node bifurcation.   Rewriting this in terms of the original 
parameters yields the two-parameter bifurcation curve in terms of $(g,I)$ :
\begin{equation}
I = I_{rh} - A_2(g)M(g)^2(g-g^*)^2 \stackrel{def}{=}I_{SN}(g) \label{ISNdef}.
\end{equation}
Thus, for fixed $g$, $I_{SN}(g)$ is the value of the current that corresponds to a saddle-node bifurcation point.  There are two things to note about $I_{SN}$.  First, since $A_2(g)$ is a strictly positive function, then $I_{SN}(g)\leq I_{rh}$ with $I_{SN}(g) = I_{rh}$ only if $g=g^*$.  Also, this curve is only defined for $g\geq g^*$. To see this, note that the saddle-node equilibrium, given by $s_{SN} = \beta = M(g)(g-g^*)$ only exists if $\beta>0$. Since $M(g)>0$, $s_{SN}$ only exists if $g\geq g^*$.  We shall see later that $g=g^*$ actually corresponds to a non-smooth co-dimension two bifurcation point.

\subsection{The Andronov-Hopf Bifurcation Condition} \label{Hopfsec}
From the analysis of subsection~\ref{stabsec}, we know that only $e_+$
may undergo a Hopf bifurcation and that det$\left.(J)\right|_{s_+}>0$
if $\beta^2+\tilde I\ne 0$.  We thus conclude that the determinant condition 
for the Hopf bifurcation is given by $I\ne I_{SN}(g)$.  To determine a 
necessary condition for the Hopf bifurcation, we begin by applying the weak 
coupling expansion \eqref{exp1} to the trace associated with $s_+$:
\begin{eqnarray}
\text{Tr}\left.(J)\right|_{s_+} 
 &=&  -\left(\frac{1}{\tau_s} + \frac{1}{\tau_w}\right) +\frac{\lambda_s}{2s_+}\left(\hat{s} g(e_r - v^*(0))-\hat{s}{g^2}{v^*}'(0)s_+-\hat{w} \right)\nonumber\\
&=& -\left(\frac{1}{\tau_s} + \frac{1}{\tau_w} + \lambda_s\hat{s}\frac{g^2}{2}{v^*}'(0)\right) +\frac{\lambda_s\hat{s}(e_r-v^*(0))}{2s_+}\left( g-\frac{\hat{w}}{\hat{s}(e_r-v^*(0))}\right) \nonumber\\
 &=&- \left(\frac{\frac{1}{\tau_s}+\frac{1}{\tau_w} + \lambda_s\hat{s}\frac{g^2}{2}{v^*}'(0)}{s_+}\right)(s_+-N(g)(g-\bar{g}))\label{trace}.
\end{eqnarray}
Note that the first term is strictly negative and $N(g)$ is a strictly positive function.  

Setting the trace to zero and using equations \eqref{tildeIdef}-\eqref{ss} which
define $\tilde I,\beta$ and $s_+$ yields
\begin{eqnarray*}
 s_+ &=& N(g)(g-\bar{g})\\
\beta + \sqrt{\beta^2 + \tilde{I}} &=& N(g)(g-\bar{g})  \\
M(g)(g-g^*) + \sqrt{M(g)^2(g-g^*)^2 + \frac{I-I_{rh}}{A_2(g)}} &=& N(g)(g-\bar{g}) \\
M(g)^2(g-g^*)^2  + \frac{I-I_{rh}}{A_2(g)} &=& (N(g)(g-\bar{g}) - M(g)(g-g^*))^2\\
\end{eqnarray*}
Solving for $I$ gives
\begin{equation}
I = I_{rh} + A_2(g)\left[N(g)^2(g-\bar{g})^2 - 2M(g)N(g)(g-\bar{g})(g-g^*) \right]
\stackrel{def}{=} I_{AH}(g)
\label{IAHdef}
\end{equation}
We thus conclude that if $I=I_{AH}(g)$ and $I\ne I_{SN}(g)$ then the equilibrium
$s_+$ has a pair of pure imaginary eigenvalues.

Recall that $N(g)$, $M(g)$, and $A_2(g)$, are positive functions. Further,
it is easy to check that $N(g)<M(g)$. This leads to several observations. 
First, since the third equation in the sequence above can only be satisfied if 
$N(g)(g-\bar g)>M(g)(g-g^*)$, it follows that if $g^*\le \bar{g}$ then no 
Hopf bifurcation occurs.  Second, from the first equation in the sequence 
above we must have $g\geq \bar{g}$ in order for the equilibrium $s_+$ to exist 
at the  Hopf bifurcation. 
When $g=\bar{g}$ $s_+=0$ and $I_{AH}=I_{rh}$.  We shall see later that the 
point $I=I_{rh}$, $g=\bar{g}$ is a codimension-2 non-smooth bifurcation point.  
Finally, if $ \bar{g} \leq g\leq g^*$, then $I_{AH}\geq I_{rh}$ with 
$I_{AH} = I_{rh}$ only if $g = \bar{g}$.  If $g>g^*$, then it is possible
for $I=I_{AH}(g)$ to intersect $I=I_{rh}$.  We denote by $\hat g$ the 
value of $g$ at the intersection point, if it exists.

We can now determine the stability of the equilibrium $e_+$ by studying
the trace equation \eqref{trace}.
Since the first term in this equation is strictly negative wherever it is 
defined (when $e_+$ exists),  the sign of the trace is determined by 
$s_+-N(g)(g-\bar{g})$.  Since $s_+$ and $N(g)$ are positive, it follows 
from the discussion above that in the case $g^*\le\bar g$ the trace negative, 
and hence $e_{+}$ is asymptotically stable, wherever $e_+$ exists. 
If $\bar g\le g^*$ then that the trace is negative (and $e_+$ is asymptotically
stable) if $g<\bar g$ or  $g>\bar g$ and $I>I_{AH}(g)$. The trace is positive 
(and $e_+$ is unstable) if $g>\bar g$ and $I<I_{AH}$.
Note that  if $I$ is sufficiently close to $I_{AH}$ then 
$e_+$ will have a pair of complex conjugate eigenvalues. 

In summary, for fixed $g$ with $g>\bar{g}$ and $\bar{g}<g^*$, the equilibrium 
$e_+$ undergoes a Hopf bifurcation at $I=I_{AH}(g)$ if $I\ne I_{SN}(g)$. We do 
not attempt to determine the criticality of this bifurcation analytically,
but will study it numerically in a later section.
Further, we can now state completely the conditions for the saddle-node
bifurcation: for fixed $g$ with $g>g^*$, the equilibria $e_+$ and $e_-$ 
undergo a saddle-node bifurcation when $I=I_{SN}(g)$ if $I\ne I_{AH}(g)$.

\subsection{The Bogdanov-Takens Bifurcation Condition}\label{BTsec} 

Recall that necessary conditions for a Bogdanov-Takens bifurcation
are $\text{det}(J)=0$ and $\text{tr}(J)=0$. Thus, from the analysis of the last
two subsections, Bogdanov-Takens bifurcations (if they exist) will occur 
at intersection point of the curves of saddle node and Hopf bifurcations 
in the $g,I$ parameter space, i.e., at values of $g$ such that
$I_{AH}(g) = I_{SN}(g)$, with $g>\max(g^*,\bar{g})$.
Using the expressions for these curves gives
\begin{eqnarray*}
I_{rh} - A_2(g)M(g)^2(g-g^*)^2 &=&  I_{rh} + A_2(g)\left[N(g)^2(g-\bar{g})^2 - 2M(g)N(g)(g-\bar{g})](g-g^*) \right]\\
0 &=& A_2(g)\left(N(g)(g-\bar{g}) - M(g)(g-g^*)\right)^2\\
&\Updownarrow&\\
0 &=& N(g)(g-\bar{g}) - M(g)(g-g^*)
\end{eqnarray*}
This latter equation may be simplified to a quadratic in $g$:
\[ a_2 g^2+a_1 g+a_0 =0\]
where
\begin{eqnarray*}
a_2&=&\frac{\hat{s}\hat{w}{v^*}'(0)}{e_r-v^*(0)}(\tau_w-\tau_s)\\
a_1&=&-\frac{2}{\tau_w}\\
a_0&=&\frac{2\eta}{\tau_s(e_r-v^*(0))}.
\end{eqnarray*}
Analysis of the quadratic equation shows that if $\tau_w<\tau_s$ or
$\tau_w>> \tau_s$ then there are no roots with $g>\max(g^*,\bar{g})$.
However, if $\tau_w>\tau_s$ and $\tau_w$ is not too large then there
can be up to two roots. Thus there may be up to two Bogdanov-Takens
points.

The case $\tau_w=\tau_s$ is degenerate. There is a single root which 
is given by $g= g^* = \bar{g}$.  This is not a standard/smooth 
Bogdanov-Takens point as the points $I_{SN}(g^*)$ and $I_{AH}(\bar{g})$ 
are not standard 
saddle-node and Hopf points.  We shall show later that this point arises
at the collision of two codimension-2 non-smooth bifurcations, and thus 
corresponds to a codimension-3 non-smooth bifurcation. 
 
\begin{figure}
\centering
\includegraphics[scale=0.5]{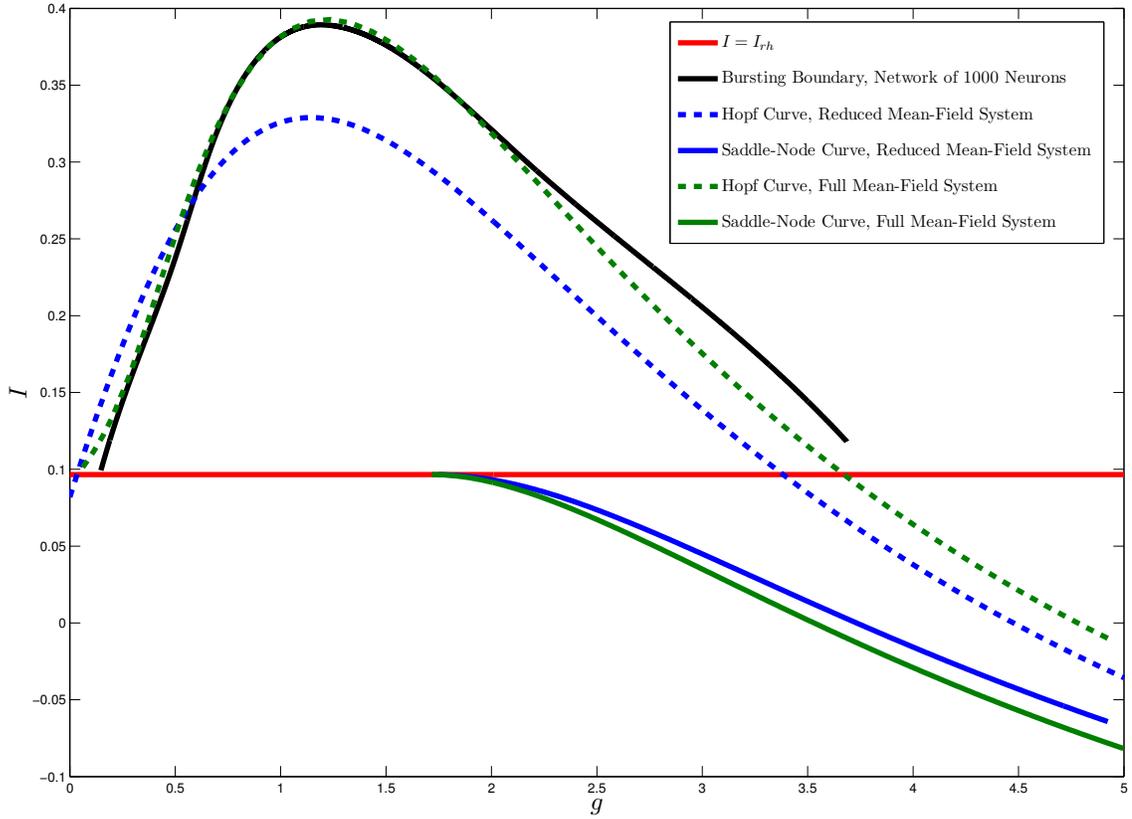}
\caption{The smooth bifurcations displayed by a network of Izhikevich neurons.  Shown in black is the bursting boundary for a network of 1000 all-to-all coupled Izhikevich neurons.  The green dashed and dotted curve are the two-parameter Hopf and saddle-node manifolds (respectively) for the full mean field equations.  These were computed using MATCONT.  The blue dashed and dotted curves are the two-parameter Hopf and saddle-node manifolds (respectively) derived from the approximate mean field equations for the Izhikevich network.  These are given by
equations \eqref{IAHdef} and \eqref{ISNdef}, respectively.  } \label{fig4}

\end{figure}

\begin{figure}
\centering
          \subfigure[Two Parameter Bifurcation Curves, Quartic Integrate and Fire mean field Models ]{\includegraphics[scale=0.5]{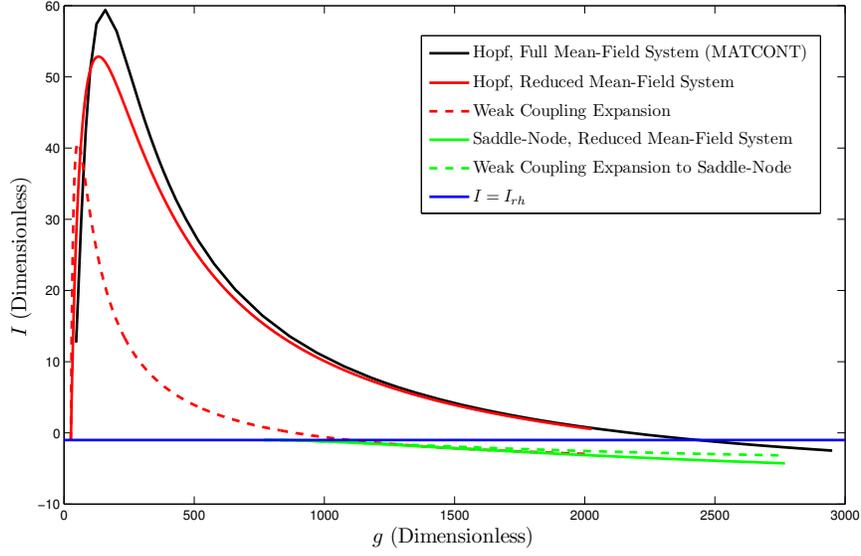}
\label{fig5a}}
\\
           \subfigure[Two Parameter Bifurcation Curves, Adaptive Exponential Integrate and Fire mean field Models ]{\includegraphics[scale=0.5]{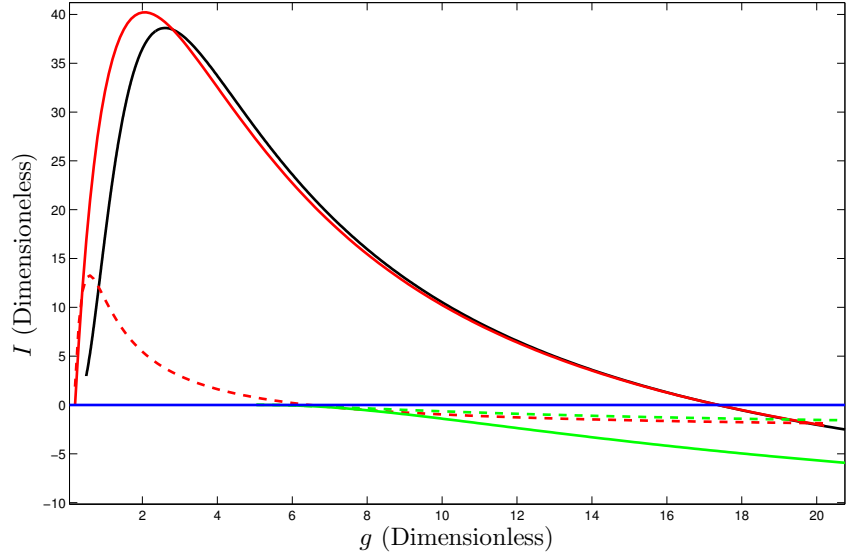}
\label{fig5b}}
        \caption{Shown above are the two-parameter Hopf bifurcation curves (in red, solid), the saddle-node bifurcation curves (green, solid), in addition to the lowest order weak coupling expansion approximation to these curves (dashed lines) for the mean field equations of the quartic integrate and fire model (a), and the Adaptive exponential integrate and fire model (b).  
The mean field system considered is the approximate system where $\langle R_i(t)\rangle \propto \sqrt{I-I^*(s,w)}$, where the value of $I^*(s,w)$ varies according to the neuron model used.   The black curve in each graph is the two-parameter Hopf bifurcation curve for the full mean field system, for comparison purposes.  This is numerically computed using MATCONT.  The bifurcation curves for the saddle-node and Hopf bifurcations are computed using the MATLAB function \emph{fsolve} on the determinant and trace equations. } \label{fig5}
\end{figure}

Figure~\ref{fig4} shows the smooth bifurcations for the mean field system
corresponding to a network of Izhikevich neurons with a the parameter values 
from \cite{us}.  Note that the Hopf-bifurcation for both the full and reduced mean-field systems corresponds closely to the onset of bursting in the actual network, as noted in \cite{us2}.  For these parameter values, $\tau_w>>\tau_s$ and no 
Bogdanov-Takens' points are observed. 
Figure~\ref{fig5} shows the smooth bifurcations for the mean field system
corresponding to networks of quartic integrate-and-fire neurons and AdEx 
neurons. In all figures the bifurcation curves derived from the weak 
coupling approximation of the model, i.e.  equations \eqref{ISNdef} 
and \eqref{IAHdef}, are compared with curves
for the full mean field model generated numerically in MATCONT \cite{MATCONT}.

\section{Non-Smooth Bifurcations} 
To study the non-smooth bifurcations for the mean field system
\eqref{MFa}--\eqref{MFc}, we will use the terminology and bifurcation 
classification for  piecewise smooth continuous systems proposed in 
\cite{nonsmooth}. 
We note that some care must be used when applying these ideas to
our system.  Letting  $x=[s,w]^T$ and recalling the definition 
\eqref{Hdef} of the switching manifold, our system may be written in the 
general form used by \cite{nonsmooth}:
\[
\dot x = \left \{ \begin{array}{ll}
F_1(x,I),& \mbox{ if } H(x,I)<0\\
F_2(x,I),& \mbox{ if } H(x,I)>0
\end{array}\right.
\]
where
\begin{eqnarray*}
 F_1(x,I)&=&\left(\begin{array}{cc}
 -\frac{s}{\tau_s}\\ -\frac{w}{\tau_w} 
\end{array}\right)\\
F_2(x,I)&=& F_1(x,I)+\sqrt{H(x,I)}\left(\begin{array}{cc}
\hat s\\ \hat w
\end{array}\right)\\
\end{eqnarray*}
However, $F_2$ is only defined for $H(x,I)>0$. 
In contrast, the work of \cite{nonsmooth} assumes that both $F_1$ and
$F_2$ are defined throughout the phase space.
Nevertheless, we able to classify a number of bifurcations in our
system by analogy with the results in \cite{nonsmooth}

We will supplement our analysis with numerical studies of our
example systems. In particular, we will perform a detailed study 
of the mean field system corresponding to a network of Izhikevich 
neurons with parameters given in Table~{\ref{table1}}.

\subsection{Boundary Equilibrium Bifurcations ($I=I_{rh}$)}
\label{BEBsec}

Recall that all the models we are considering 
have an equilibrium $e_0=(0,0)$ which exists (and is a stable node)
if $I<I^*(0,0)=I_{rh}$.  When $I=I_{rh}$ this equilibrium lies on
the switching manifold $I-I^*(s,w)=0$. When  $I>I_{rh}$,  this equilibrium
no longer exists as the origin is not an equilibrium of part of the 
mean field system corresponding to $I-I^*(w,s)>0$.  In the terminology of 
non-smooth systems, the origin is a {\em virtual equilibrium} 
of the system for $I>I_{rh}$  and undergoes a {\em boundary equilibrium 
bifurcation} (BEB) when $I=I_{rh}$. The exact nature of this bifurcation
depends on the value of $g$, in particular, its relationship to $g^*$, 
$\bar g$ and $\hat g$.

To determine the nature of the boundary equilibrium bifurcation,
we begin by studying the nontrivial equilibria $e_{\pm}=(s_\pm,\eta s_\pm)$ 
when $I=I_{rh}$. Recalling the form \eqref{ss} for $s_\pm$ and noting
that $\tilde I=0$ when $I=I_{rh}$, we find
\begin{eqnarray*}
s_+(\beta,0) &=&\begin{cases} 0 & \beta<0 \\ 2\beta & \beta \geq 0 \end{cases} \\
s_-(\beta,0) &=& \begin{cases} 2\beta & \beta <0 \\ 0 & \beta \geq 0  \end{cases}
\end{eqnarray*}
Thus for $g<g^*$,  $e_+$ collides with $e_0$ at $I=I_{rh}$, and for 
$g>g^*$, $e_-$ collides with $e_0$. 

Consider first the case $g^*<\bar g$ (which corresponds to $\tau_w<\tau_s$).
In this case there is no Hopf bifurcation, so the results are straight forward.
When $g<g^*$ $e_+$ is a sink which exists for $I>I_{rh}$. 
It collides with $e_0$ when $I=I_{rh}$ and ceases to exist
when $I<I_{rh}$.  Putting this together with the description of the
existence and stability results for $e_0$, we conclude that, for this 
range of $g$ values, the system undergoes a {\bf persistence} BEB
at $I=I_{rh}$. This will be either a focus/node or node/node persistence
BEB depending on the classification of $e_+$. 
When $g>g^*$, recall that the equilibrium $e_-$ 
is a saddle when it exists (for $I_{SN}<I<I_{rh}$).  Since $e_0$ also 
exists for $I<I_{rh}$ and is a stable node, we conclude that for $g>g*$ 
there is a non-smooth saddle node BEB at $I=I_{rh}$.  
 
\begin{figure}
\centering

          \subfigure[$g<\bar{g}$ ]{\includegraphics[scale=0.53]{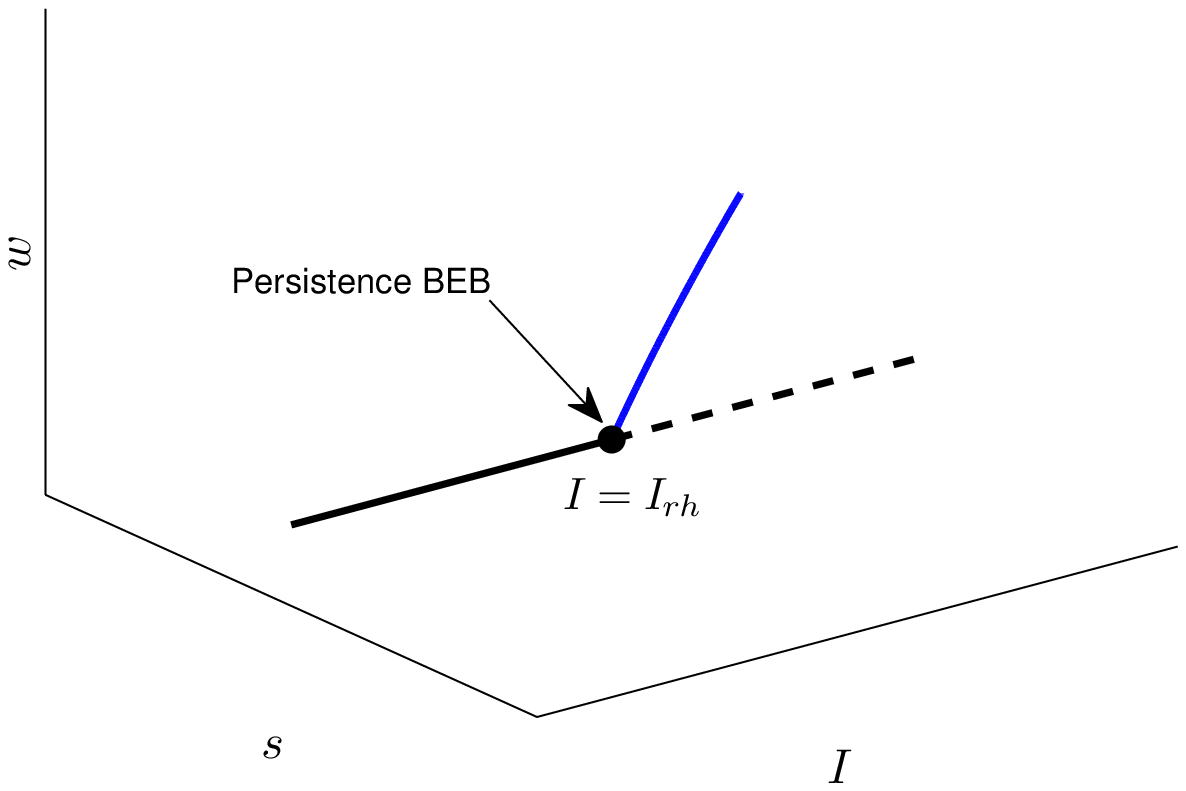}
\label{BEBa}}
\quad
           \subfigure[ $\bar{g}<g< g^*$]{\includegraphics[scale=0.53]{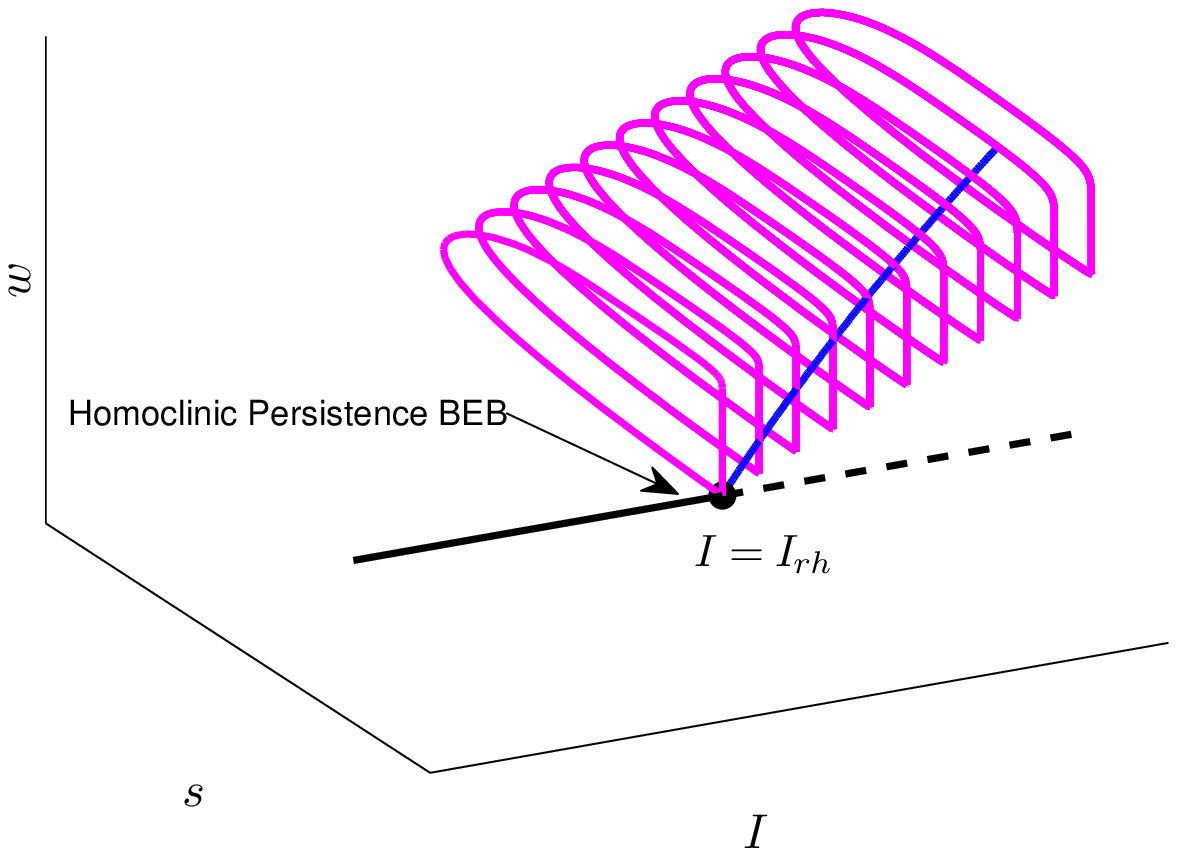}
\label{BEBb}}
\\
          \subfigure[ $g^*<g<\hat g$  ]{\includegraphics[scale=0.53]{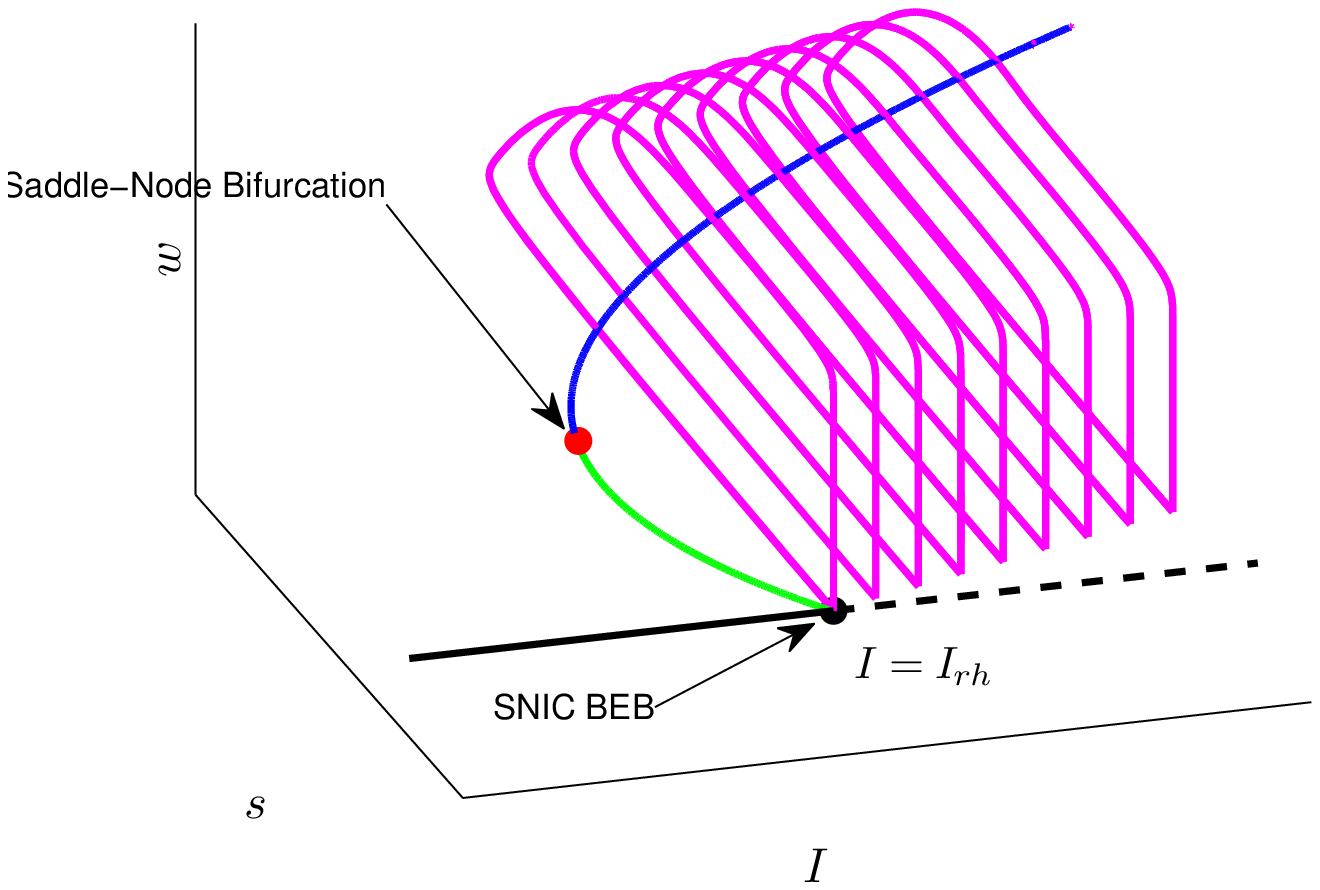}
\label{BEBc}}
\quad
           \subfigure[ $g\gg \hat g$ ]{\includegraphics[scale=0.53]{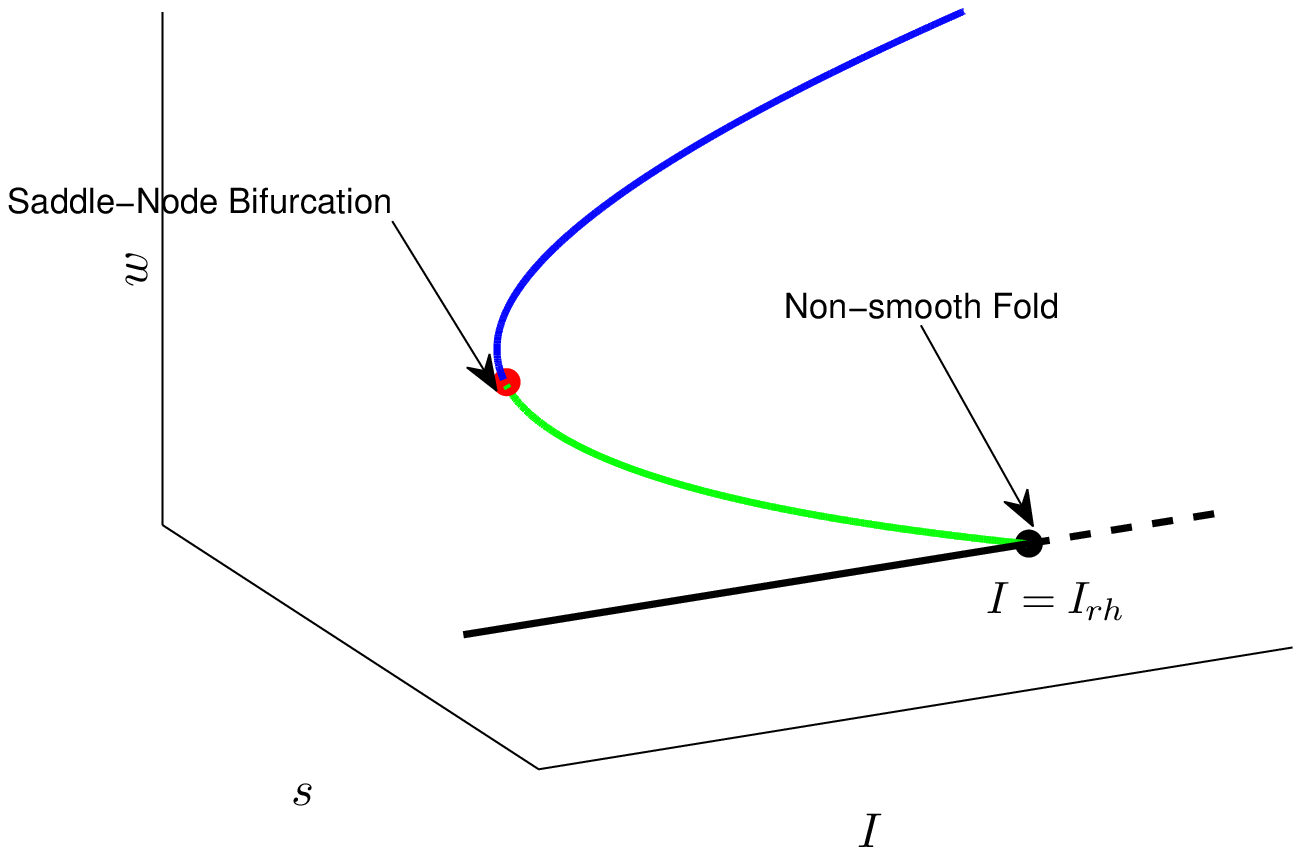}
\label{BEBd}}
\\
    \caption{Shown above are the four branches of boundary equilibrium bifurcations (BEB) that have been found in the mean field system for the Izhikevich network.  In all figures, the equilibria are $e_0$ (black), $e_+$ (blue) and $e_-$ (green), and solid lines indicate real equilibria, while dashed lines indicate virtual ones.  The magenta lines are the non-smooth limit cycles determined via direct numerical integration.  (a) the equilibrium $e_+$ collides with $e_0$ at $I=\frac{\alpha^2}{4}$.  This results in the disappearance of $e_+$ for $I<\frac{\alpha^2}{4}$, while $e_0$ persists as a stable node.  The case is identical for (b), except however that the non-smooth limit cycle collides with the BEB equilibrium point in a kind of non-smooth homoclinic bifurcation.  (c) the equilibrium $e_-$ exists and is an unstable saddle for $I<\frac{\alpha^2}{4}$, as does the stable node $e_0$.  These equilibria collide in a boundary equilibrium bifurcation at $I=\frac{\alpha^2}{4}$, and $e_-$ is destroyed while $e_0$ becomes virtual.  The bifurcation diagram in (d) is similar to that in  (c) except for the emergence of a homoclinic limit cycle at the bifurcation point in a kind of non-smooth SNIC bifurcation.} \label{figBEB}
\end{figure}

Now consider the case $\bar{g}<g^*$.
For $g<\bar{g}$, analysis similar to that above shows 
the system undergoes a {\bf persistence BEB} at $I=I_{rh}$.
Figure~\ref{BEBa} shows this bifurcation for the mean field system 
corresponding to the Izhikevich network with parameters as in Table \ref{table1}.

The situation for $\bar g<g<g^*$ is similar, except that $e_+$ is now an
unstable focus for $I>I_{rh}$.  Thus for this range of $g$ values, there
is a focus/node persistence BEB at $I=I_{rh}$. Since $e_+$ is a source
and $e_0$ is a sink, we may expect (by analogy with the results in 
\cite{nonsmooth}) that a stable non-smooth limit cycle surrounding $e_+$ 
will be created as $I$ increases through $I_{rh}$. Figure~\ref{BEBb} confirms
this for the mean field system corresponding to the 
Izhikevich network. Note that in this example, the amplitude of
the limit cycle does not got to zero as $I$ approaches $I_{rh}$.
(See also in Figure~\ref{fig8c})  Further, the period
of the limit cycle diverges as $I\rightarrow I_{rh}^+$. See
Figure~\ref{fig8a}.  Thus the limit cycle appears to be created 
in homoclinic-like bifurcation as $I$ increases through $I_{rh}$.
We will thus refer to this as a {\bf homoclinic persistence BEB}.

When $g>g^*$, analysis similar to that above shows that 
there is a non-smooth saddle node BEB at $I=I_{rh}$.  Based on the 
analysis of the equilibrium
points, there is no reason to expect anything more to occur with 
this bifurcation.  However, our numerical examples show two cases. 
Figure~\ref{BEBd} shows that a simple non-smooth saddle-node BEB 
occurs for the mean field system corresponding to the Izhikevich 
network with $g\gg \hat g$.  Figure~\ref{BEBc} shows the bifurcation 
for the same system with $g^*<g<\hat{g}$.  In this case 
there is a non-smooth limit cycle for $I>I_{rh}$ that appears to 
be destroyed when $I=I_{rh}$.  Thus this bifurcation appears
to be a non-smooth version of the Saddle-node on an invariant
circle (SNIC) bifurcation. We will refer to it as a {\bf SNIC BEB}.
The transition between the two types of BEBs that occur for
$g>g^*$ will be discussed in a later section.

Based on our numerical results we hypothesize that a 
non-smooth limit cycle may be destroyed in a homoclinic-like 
bifurcation as $I$ decreases through $I_{rh}$.  We support
this hypothesis in two ways.

First, consider the vector field in
the neighbourhood of the origin.  Recall that the origin is always 
an attractor when it lies in the region where $H(s,w)<0$.  
In the region where $H(s,w)>0$, setting $I=I_{rh}$ and retaining only
the highest order terms in $s$ and $w$ gives:
\begin{eqnarray*}
s' &=& -\frac{s}{\tau_s} + \hat{s} \sqrt{gs(e_r-v^*(0)) - {v^*}'(0)(gs)^2/2 -w }  \approx \hat{s}\sqrt{gs(e_r-v^*(0))-w} \\
w' &=& -\frac{w}{\tau_w} + \hat{s} \sqrt{gs(e_r-v^*(0)) - {v^*}'(0)(gs)^2/2 -w } \approx \hat{w}\sqrt{gs(e_r-v^*(0))-w}
\end{eqnarray*}
Thus, for $0<s,w\ll 1$, and $I>I^*(s,w)$ the vector field point away from 
the origin and the boundary equilibrium $(0,0)$ is a repeller in this region.
Since the boundary equilibrium point is as a repeller on one side of the switching manifold and and an attractor on the other, it is possible for a non-smooth
homoclinic orbit to this equilibrium point to exist when $I=I_{rh}$. 

Second, we show that under certain parameter conditions, if a non-smooth 
limit cycle surrounds the equilibrium $e_{+}$, it must be destroyed 
when $I=I_{rh}$. To do this we show that trajectories that cross 
the switching manifold when $I_{rh}$ lie within the basin of attraction of the origin.  Thus any non-smooth limit cycle must become homoclinic to the 
origin at $I=I_{rh}$.  Note that if  $I-I^*(s,w)<0$, then 
\begin{eqnarray*}
\frac{dw}{ds} = \frac{\tau_s}{\tau_w} \frac{w}{s} =\gamma \frac{w}{s}
\end{eqnarray*}
and thus $w = Cs^\gamma$ for some constant $C$.  Assuming that the trajectory starts with $(s_0,w_0)$ on the switching manifold then 
$w = w_0 \left(\frac{s}{s_0}\right)^\gamma$ where 
$w_0 = gs_0(e_r-v^*(0))-{v^*}'(0)\frac{(gs)^2}{2}$.
Now suppose this trajectory crosses the switching manifold again at $(s,w)$.
Then
\[
 w_0 \left(\frac{s}{s_0}\right)^\gamma = gs(e_r-v^*(0)) - {v^*}'(0)\frac{(gs)^2}{2} 
\]
Clearly two solutions of this equation are $(s_0,w_0)$ and $(0,0)$.   
Dividing through by $s$ and simplifying one obtains 
\begin{equation}
(1-ks_0)\frac{s^{\gamma-1}}{s_0^{\gamma-1}} = 1-ks
\label{homeq}
\end{equation}
where $k=\frac{g{v^*}'(0)}{2(e_r-v^*(0))}$.

If $\gamma>1$ the left hand side of \eqref{homeq} is monotonically 
increasing while the right
hand side is a line with negative slope. Hence $(s_0,w_0)$ is the unique 
intersection point. This means every trajectory that enters the region 
$I-I^*(s,w)<0$ when $I=I_{rh}$ is attracted to the origin.
If $\gamma < 1$, the left side of \eqref{homeq} is now
monotonically decreasing. Unless the line is tangent the curve at $(s_0,1)$
there will always be another intersection point. 
Rearranging the equation shows that this intersection point
will occur for $s<s_0$ if $gs_0$ is sufficiently small. For fixed $g$, 
this means that any trajectory that starts on the switching manifold at 
$(w_0,s_0)$ with $s_0$ sufficiently small will be attracted to the origin.
Thus all non-smooth limit cycles that are close enough to the origin 
for $I>I_{rh}$ will become homoclinic to the origin when $I=I_{rh}$.



Given how $g=\bar{g}$, and $g=g^*$ delimit the different types of BEB bifurcations, it should be clear that these special points represent higher codimension bifurcations along the $I=I_{rh}$ line.  We shall explore these bifurcations further below, in addition to determining the geometrical meaning of these points.

\subsection{Saddle-Node Boundary Equilibrium Bifurcation ($I=I_{rh},g=g^*$)} 

From the previous section, we concluded that the point $I=I_{rh}, g=g^*$ is 
a special codimension-two bifurcation point where the boundary equilibrium bifurcation (BEB) changes from a persistence BEB to a non-smooth
saddle-node.  Note that the smooth branch of saddle-node bifurcations found earlier actually emanates out from the codimension-2 point $(g^*,I_{rh})$.  
We will show here that it does so in a highly non-generic way as the 
saddle-node equilibrium hits switching manifold tangentially at the BEB, and 
is the only equilibrium point that interacts with the switching manifold in this way.  

We have seen that regardless of the parameter values, all the nontrivial 
equilibria lie on the curve $w = \eta s$.  Thus as any parameter is varied
the nontrivial equilibrium will follow this curve, which has slope
\begin{equation}
 w'(s)=\eta.
\label{eqmslope}
\end{equation}
Further, the only equilibrium that can be boundary equilibrium point is 
$e_0 = (0,0)$, the non-firing solution.  Now the switching manifold can be 
written as 
\[ w(s) = I+F(v^*(s))+gs(e_r-v^*(s)) \]
Thus, at the slope of the switching manifold at the BEB is
\begin{equation}
w'(0) = g(e_r - v^*(0))
\label{SMslope}
\end{equation} 
Equating \eqref{eqmslope} and \eqref{SMslope} shows that the nontrivial 
equilibrium undergoing the BEB will hit the switching manifold tangentially
only if $g = g^* =  \frac{\eta}{e_r - v^*(0)}$. 
From this it is straightforward
to show that with $g=g^*$ fixed, the nontrivial equilibrium $e_+$ hits the
switching manifold tangentially as $I\rightarrow I_{rh}$ $s_+\rightarrow 0$. 
More interesting is to consider what happens when $g$ is varied. 
From our previous analysis we know that at the saddle-node bifurcation point, 
the saddle-node equilibrium, $e_{SN}=(s_SN(g),\eta s_SN(g))$ is defined by 
\begin{eqnarray*}
s_{SN}(g) = M(g)(g-g^*)
\end{eqnarray*}
Thus, as $g\rightarrow g^*$, $s_{SN}(g)\rightarrow 0$.   This implies that 
the saddle-node equilibrium hits the switching manifold tangentially 
at $g=g^*$, $I=I_{rh}$. 

In summary the point $g = g^*, I= I_{rh}$ is the collision between three 
branches of codimension-1 bifurcations: a pair of non-smooth boundary
equilibrium bifurcations and a smooth branch of saddle-node bifurcations. 
The details of the BEB involve depend on the relationship between $g^*$
and $\bar g$.  If $g^*<\bar{g}$, BEBs are simple: a simple node/focus 
or focus/focus persistence BEB occurs for $g<g^*$ and a non-smooth 
saddle-node BEB occurs for $g<g^*$. The case $\bar{g}<g^*$ is more
complex due to the possible presence of limit cycles associated with
the Hopf bifurcation. In the case we studied numerically and described 
in section~\ref{BEBsec}, for $g<g^*$ 
we observe a homoclinic persistence BEB and for $g>g^*$ we observe a SNIC BEB. 

While this bifurcation may be complicated, the bifurcation point can be
determined analytically for all the models.  Of particular interest is the 
underlying physical interpretation.  Associated with this point
is a region in the $g>g^*$, $I<I_{rh}$ quadrant of the parameter space with both 
stable firing and non-firing solutions, and hence bistability.  
This bifurcation point is shown in detail in figure \ref{fig10b}. 
 
\subsection{Limit Cycle Grazing Bifurcation} 

The Andronov-Hopf bifurcation described in section~\ref{Hopfsec} leads to the
creation of a limit cycle.  As $I$ moves away from the bifurcation point,
the amplitude of the limit cycle may increase enough that it hits 
the switching manifold tangentially, resulting in a grazing bifurcation.
It is difficult to say much in general about the nature of this 
bifurcation, however, analysis similar to that in the last section
shows that if $I<I_{rh}$ then once a trajectory enters the region
$I-I^*(s,w)<0$, it cannot leave, but will be attracted to the origin.
Thus we expect that if a grazing bifurcation occurs for $I<I_{rh}$
it will lead to the destruction of the limit cycle.

\begin{figure}
\centering
          \subfigure[Grazing Bifurcation, Persistence, and Non-Smooth Saddle-Node of Limit Cycles]{\includegraphics[scale=0.45]{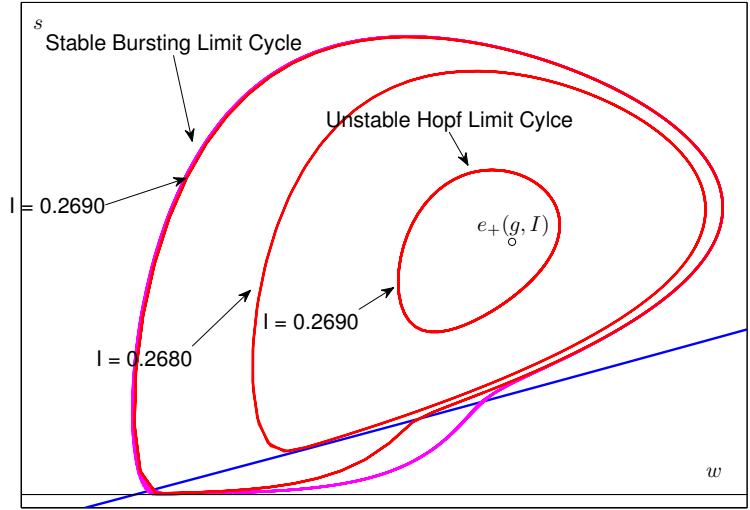}
\label{fig7a}}\\
          \subfigure[Grazing Bifurcation, Destruction ]{\includegraphics[scale=0.45]{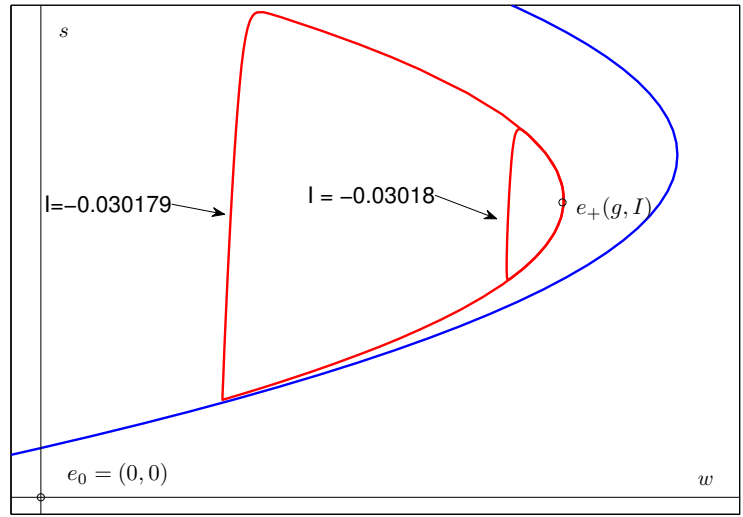}
\label{fig7b}}

\caption{Limit cycle grazing bifurcations for the Izhikevich system.  (a) As we increase $I$ above $I_{AH}(g)$, for fixed $g$, the unstable limit cycle (shown in red) generated by the sub-critical Hopf bifurcation increases in amplitude.  For large enough $I$, the limit cycle grazes the switching manifold (shown in blue).  After the grazing, the limit cycle becomes non-smooth and subsequently collides with the non-smooth stable limit cycle (shown in pink).  The two limit cycles annihilate each other in a non-smooth saddle node of limit cycles.  Note that as we vary $I$, the switching manifold, the point $e_+$, and the unstable limit cycle all vary.  However, aside from the unstable limit cycle, these other sets do not vary significantly.  Thus, for clarity, we have only shown the switching manifold and stable non-smooth limit cycle for $I=0.2690$, and $e_+$ for $I=0.2604$. 
(b) For $I<I_{rh}$ the grazing bifurcation destroys the limit cycle.}\label{fig7}
\end{figure}

\begin{figure}
\centering
          \subfigure[Amplitude of the Stable Non-smooth Limit Cycle]{\includegraphics[scale=0.5]{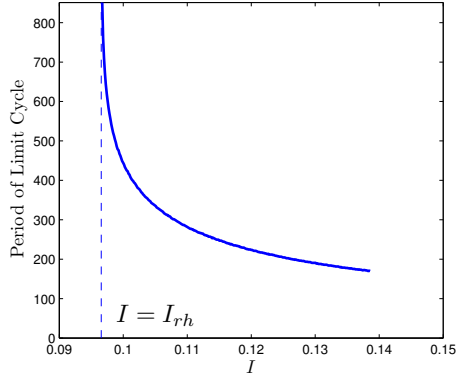}
\label{fig8a}}\quad
\centering
          \subfigure[Period of the Stable Non-smooth Limit Cycle ]{\includegraphics[scale=0.5]{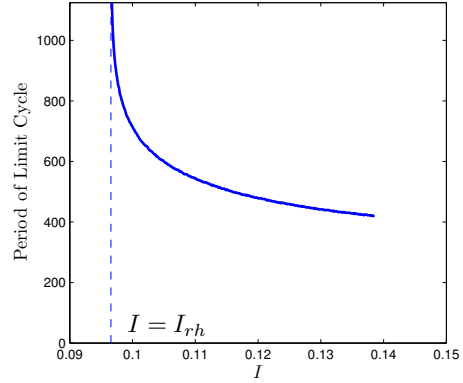}
\label{fig8b}}
\\
          \subfigure[Amplitude of the Stable Non-smooth Limit Cycle]{\includegraphics[scale=0.5]{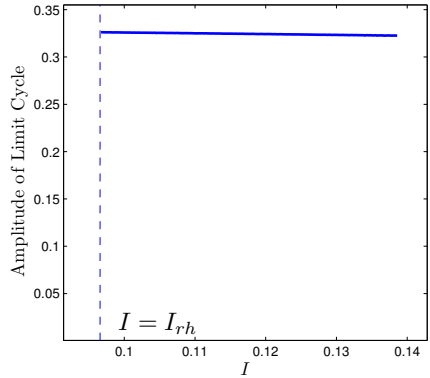}
\label{fig8c}}\quad
          \subfigure[Period of the Stable Non-smooth Limit Cycle ]{\includegraphics[scale=0.5]{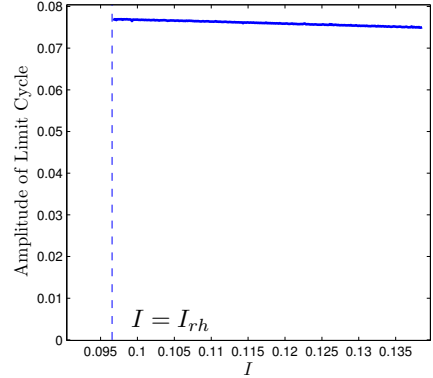}
\label{fig8d}}

\caption{Shown above is the amplitude (a) and period (b) of the bursting limit cycle for fixed $g$ with $\bar{g}<g<g^*$ (left column) and $g>g^*$ (right column) as $I\rightarrow I_{rh}$.  These two quantities are resolved via direct numerical simulation of the limit cycle.  Note the period diverges as $I\rightarrow I_{rh}$, while the amplitude is non-zero, indicative of a homoclinic limit cycle.  The amplitude is computed as the difference between the maximum and minimum $w$ component in the steady state limit cycle.    }\label{fig8}
\end{figure}

To gain more insight, we performed a numerical study of the
mean field system corresponding to the Izhikevich network with 
parameter values as in Table \ref{table1}. We first confirmed that the 
Hopf bifurcation is subcritical, using MATCONT and by numerically 
simulating the time reversed system.  We then showed
that the unstable limit cycle generated  by the Hopf can undergo
two different types of grazing bifurcations,
depending on the value $I$. 
For $I>I_{rh}$, the grazing bifurcation that occurs 
is a persistence type grazing, i.e., the unstable limit cycle generated 
via the subcritical Hopf bifurcation just becomes non-smooth after
the grazing bifurcation.  This is shown in Figure \ref{fig7a}.  Here, the limit cycle undergoes a grazing bifurcation at $I=0.2680$, and it persists past it.  Its amplitude rapidly increases past the grazing bifurcation, 
and it almost immediately undergoes a non-smooth saddle-node of limit cycles 
with a stable non-smooth limit cycle. 
For $I<I_{rh}$, the grazing bifurcation is a destruction type grazing as
the limit cycle ceases to exist after the grazing for the reason discussed 
above.  This is shown in Figure \ref{fig7b}.  

If the Hopf were supercritical we would expect to see the same two 
types of grazing bifurcations. The only difference would be that 
the grazing bifurcation would occur for $I<I_{AH}$ and we would
not expect the saddle-node of limit cycles bifurcation to occur.

\subsection{Hopf Boundary Equilibrium Bifurcation ($I=I_{rh},g=\bar{g}$)} 

The analysis of section~\ref{BEBsec} showed that when $\bar{g}<g^*$ the 
point $I=I_{rh},g=\bar{g}$
is a codimension-two bifurcation point where the boundary equilibrium changes
from a simple focus/node persistence BEB to a homoclinic persistence BEB.
Recall that the two parameter Hopf bifurcation curve is given by $I=I_{AH}(g)$
as defined in section~\ref{Hopfsec}. From the analysis in that section,
the equilibrium point on the Hopf curve is $e_{AH}=(s_{AH},\eta s_{AH})$ 
where $s_{AH}(g) = N(g)(g-\bar{g})$.  Setting $I=I_{AH}(g)$ we see that 
as $g\rightarrow\bar{g}$, $I_{AH}\rightarrow I_{rh}$ and 
$e_{AH}\rightarrow e_0$, that is the Hopf equilibrium point undergoes a BEB
at $I=I_{rh},g=\bar{g}$.  We thus refer to this point as a Hopf boundary 
equilibrium bifurcation (Hopf BEB).

An alternative way to characterize the Hopf BEB is to fix $I=I_{rh}$ and
let $g\rightarrow \bar{g}^+$.  On $I=I_{rh}$, the mean field system for the 
Izhikevich network may be approximated as follows: 
\begin{eqnarray*}
s' &=& -\frac{s}{\tau_s} + \hat{s} \sqrt{gs(e_r-\frac{\alpha}{2})-\frac{(gs)^2}{4} -w} \approx \hat{s} \sqrt{gs(e_r-\frac{\alpha}{2})-w}\\
w' &=& -\frac{w}{\tau_w} + \hat{w} \sqrt{gs(e_r-\frac{\alpha}{2})-\frac{(gs)^2}{4} -w}\approx \hat{w} \sqrt{gs(e_r-\frac{\alpha}{2})-w}\\
\end{eqnarray*}
for $(s,w)$ in the vicinity of the origin.  Thus, we have 
\begin{eqnarray*}
\frac{dw}{ds} &=& \frac{\hat{w}}{\hat{s}}+ H.O.T.\\
\Rightarrow w &=&\frac{\hat{w}}{\hat{s}}s + H.O.T.
\end{eqnarray*}
for the trajectory of the homoclinic limit cycle.   Additionally, linearizing the switching manifold about the origin yields; 
\begin{eqnarray*}
w = gs\left(e_r-\frac{\alpha}{2}\right)
\end{eqnarray*}
Now, using these two equations we can solve for grazing bifurcations of the homoclinic limit cycle with the switching manifold at the origin.  Solving the grazing condition $w'(0) = \frac{\hat{w}}{\hat{s}}$ yields
\begin{eqnarray*}
g = \frac{\hat{w}}{\hat{s}(e_r-\alpha/2)}= \bar{g}
\end{eqnarray*} 
Thus, the Hopf BEB bifurcation can be seen as a grazing bifurcation which destroys the non-smooth homoclinic limit cycle to the origin.  

Our analysis so far shows three branches of bifurcation emanating from 
this codimension-two point: two non-smooth BEB branches and a branch of 
Hopf bifurcation.  As shown in Figure~\ref{BEBfig}, for $g<\bar{g}$ there is a simple persistence BEB, while for $\bar{g}<g<g^*$ 
there is a homoclinic persistence BEB. 
We have numerically studied the bifurcations that occur in 
a neighbourhood of this point for the Izhikevich model 
and find that that two more branches of bifurcation appear to 
emanate from this point we describe below. 

Let $g$ be fixed with $g>\bar{g}$ and consider the sequence of bifurcations 
involving limit cycles.  At $I=I_{rh}$ a stable non-smooth limit cycle is 
created in a homoclinic persistence BEB, at $I=I_{AH}>I_{rh}$ an unstable 
smooth limit cycle is created in a subcritical Hopf bifurcation. As
$I$ increases the smooth limit cycle becomes non-smooth in a grazing
bifurcation and then is destroyed along with the stable non-smooth limit
cycle in a saddle-node of limit cycles. We wish to determine how the
grazing and saddle-node of limit cycles bifurcations behave near $g=\bar{g}$.

\begin{figure}
\centering
          \subfigure[Amplitude of the Stable Non-smooth Limit Cycle ]{\includegraphics[scale=0.7]{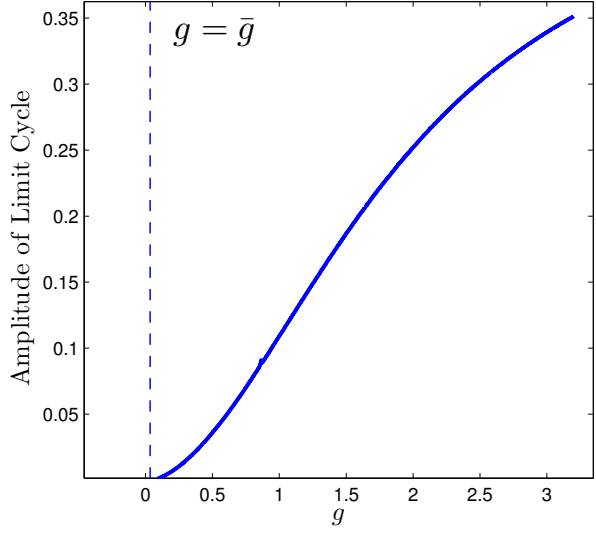}
\label{fig9a}}\\
\centering
          \subfigure[Period of the Stable Non-smooth Limit Cycle ]{\includegraphics[scale=0.73]{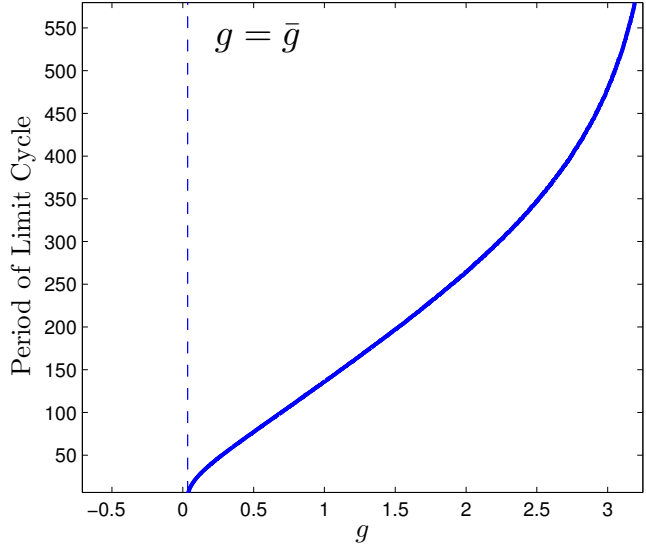}
\label{fig9b}}
\caption{Shown above is the amplitude (a) and period (b) of the bursting limit cycle followed along the two-parameter Hopf bifurcation curve.  
The Hopf bifurcation curve is entirely parameterized by $g$, in the $(I,g)$ plane, and thus as we decrease $g$, we can compute the amplitude and period of the bursting limit cycle via direct numerical simulations.  As can be seen, the amplitude decreases towards 0 as $g\rightarrow \frac{\hat{w}}{\hat{s}(e_r-\alpha/2) } = \bar{g}$, as does the period.  As the bursting limit cycle is the exterior limit cycle in a non-smooth saddle node bifurcation of limit cycles, this bifurcation must also emanate from $\bar{g}$.  Additionally, as the saddle-node of limit cycles occurs subsequent to a persistent grazing bifurcation of the unstable Hopf limit cycle, the grazing bifurcation must also emerge from this point.   Also note that this is the only point in the parameter space where the homoclinic limit cycle generated does not have a divergent period as $I\rightarrow I_{rh}$.  This is due to the fact that the homoclinic limit cycle has collapsed down to a point exactly at $g=\bar{g}$, and thus does not exist at this parameter value.  }\label{fig9}

\end{figure}

\begin{figure}
\centering
          \subfigure[Hopf BEB Point (Bottom Left Corner)) ]{\includegraphics[scale=0.53]{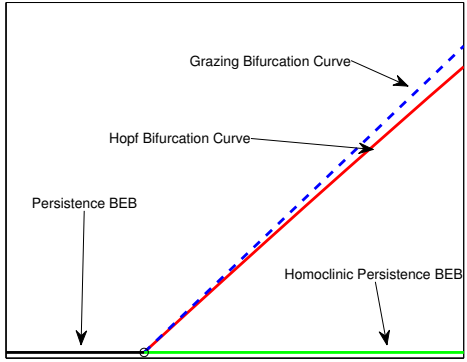}
\label{fig10a}}
\quad
           \subfigure[Saddle-Node BEB (Center) ]{\includegraphics[scale=0.53]{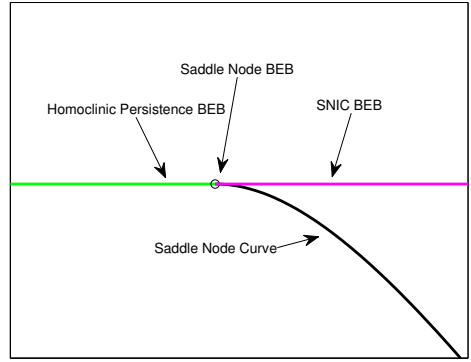}
\label{fig10b}}
\\
          \subfigure[Grazing Alternation (Bottom Right Corner) ]{\includegraphics[scale=0.53]{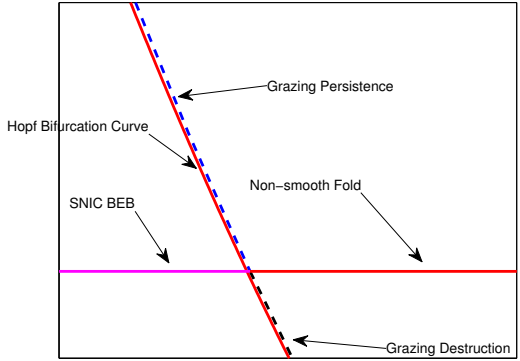}
\label{fig10c}}
\quad
           \subfigure[Total Bifurcation Diagram]{\includegraphics[scale=0.53]{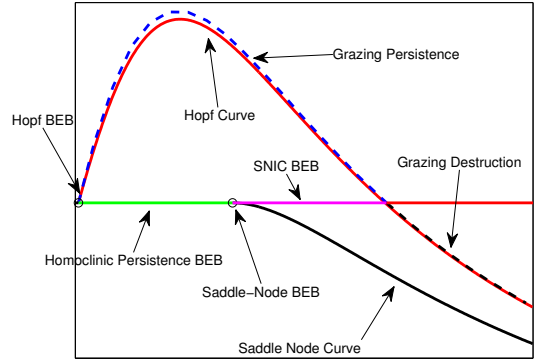}
\label{fig10d}}
\\

        \caption{Shown above is the entire bifurcation sequence for the Izhikevich model, including all known non-smooth and smooth bifurcation points.  (d) is the entire diagram in the two-parameter space. (a), (b), and (c) are the bottom left, center, and bottom right regions, respectively.   
(a) The co-dimension two bifurcation point involving the collision of a branch of Hopf bifurcations with the switching manifold.  This co-dimension two point also involves a collision with a branch of grazing bifurcations of the unstable limit cycle generated by the sub-critical Hopf, in addition to a branch of saddle-node of limit cycles (not shown for clarity).  A non-smooth SNIC bifurcation, and BEB persistence bifurcation also collide simultaneously at the codimension two point $(\frac{\hat{w}}{\hat{s}(e_r-\alpha/2)},\frac{\alpha^2}{4})$.  
(b) The codimension two saddle-node grazing point, which occurs when a saddle-node bifurcation grazes a switching manifold.  The saddle-node branch of bifurcations collides at the codimension-two point $(\frac{\eta}{e_r-\alpha/2},\frac{\alpha^2}{4})$ along with two branches of non-smooth SNIC bifurcations.  
(c) A global codimension-two point.  This bifurcation point involves the switching of a grazing bifurcation in the unstable Hopf limit cycle from a persistence case, to a destruction case.  The non-smooth SNIC bifurcation also collides with a branch of BEB persistence bifurcations for the equilibrium $e_-(g,I)$.
}    \label{fig10}
\end{figure}

To do this we followed the stable non-smooth limit cycle along the Hopf
bifurcation curve. Specifically, we numerically computed the amplitude and
period of the limit cycle along the curve $(g,I_{AH}(g))$ in the $(g,I)$
parameter space with $g\rightarrow \bar{g}$.
The results are shown in Figure \ref{fig9}, specifically figure \ref{fig9a}.  The stable non-smooth limit cycle is computed using direct simulations of the ODE system, where the system is initialized exterior to the limit cycle in the phase plane which ensures convergence.  From this figure, we can see that the 
amplitude of the stable non-smooth limit cycles goes to 0 as 
$g\rightarrow \bar{g}$.  This implies that this limit cycle collapses to the origin $(0,0)$.   But as this bursting limit cycle is one part of the saddle-node of limit cycles bifurcation, then this bifurcation must also emerge from 
Hopf BEB.  
Since the grazing bifurcation lies between the saddle-node of limit cycles
and the Hopf bifurcation,  the persistence grazing bifurcation must also emerge from the point $g=\bar{g},I=I_{rh}$.  
The entire sequence of bifurcations near the Hopf BEB is shown in figures \ref{fig10a} and \ref{fig10d}.   


\subsection{A Co-dimension 3 Non-smooth Bifurcation} 

We briefly note that if $\tau_w = \tau_s$, then we have 
\begin{equation}
\bar{g} = g^*,
\end{equation}
which means that the Hopf and saddle-node BEB points coincide in a non-smooth codimension-3 bifurcation point.   This bifurcation point may be thought of as a Bogdanov-Takens equilibrium point lying on a switching manifold. 
However, we note that there is no Bogdanov-Takens bifurcation (or for that matter saddle-node or Hopf bifurcations) at this point in the classical sense, as the Jacobian of the system diverges, and hence the conditions associated with these different smooth bifurcations cannot be satisfied. 

This point appears to act as an organizing center for the bifurcation diagram, with
all the non-smooth bifurcations emanating from it.  Due to the complexity of this point, 
we will leave its analysis for future work.  However, it does illustrate how rich the non-smooth bifurcation sequence of this relatively simple PWSC system is.  

\subsection{A Global Co-dimension 2 Non-smooth Bifurcation} 

In addition to the two local non-smooth bifurcations that occur at $g=\bar{g}$, and $g=g^*$, there appears to be a global codimension-2 bifurcation that occurs for these 
mean field systems.  Recall that there are two different types of grazing bifurcations, 
a destruction type (which occurs for $I<I_{rh}$) and a persistence type (which occurs
for $I>I_{rh}$. These are shown in figure \ref{fig7}. Thus there is a co-dimension
two point when the grazing bifurcation crosses $I=I_{rh}$.  As for the other
co-dimension two points, one may expect that would be a change in the BEB bifurcations 
at this point. In the case we have studied numerically it appears
that the BEB changes from SNIC type before this transition to a regular non-smooth
fold after.  This is shown in figure \ref{fig10c} and figure \ref{fig10d}. Note that
this transition occurs for $g>\hat{g}$, i.e., after the second intersection of the
Hopf curve with $I=I_{rh}$.
It also appears that the saddle-node of non-smooth limit cycles bifurcation emanates 
from this point.
Note that this does not imply that there is a second impact with the Hopf equilibrium and the switching manifold, as $s_{AH}(g) = N(g)(g-\bar{g})>0$.  This bifurcation results in the destruction of the homoclinic limit cycle that exists on $I=I_{rh}$, and it is very difficult to analyze, as it is a non-local co-dimension 2 non-smooth bifurcation.   Geometrically, however it occurs when the unstable smooth limit cycle (generated via the Hopf bifurcation) grazes the switching manifold at $I=I_{rh}$. If the Hopf bifurcation were supercritical 
instead of subcritical we would expect a similar codimension two point to occur (if a
grazing bifurcation occurred). However, it would occur for $g<\hat g$.

Again, due to the complexity of this particular bifurcation, further analysis is beyond the scope of this paper, and we leave it for future work.  

\section{Non-Smooth Bifurcations Demonstrated in the Network Simulations} 

While the preceding analysis revealed a great deal of novelty and non-smooth bifurcations for the reduced mean-field system, in order for the non-smooth analysis to be useful, it has to be reflective of the phenomenon displayed by the actual network.  Here, we will demonstrate many of the non-smooth bifurcations predicted in the analysis are present in a full network of neurons. 

Unfortunately however, one cannot easily simulate the large network of neurons in such a way as to expose unstable equilibria and limit cycles.  For example, the equilibrium point $e_+$ is a saddle in the mean-field, and short of somehow initializing the network of neurons on the stable manifold of the saddle, it cannot be resolved via direct simulations.  However, the unstable node $e_+$ can be resolved by modifying the network as follows. Using the separation of time scales between the fast variable $s$, and the slow variable $w$, we replace the full network 
\eqref{v_fullnet}--\eqref{wjump_fullnet} by the following:
\begin{eqnarray*}
\dot{v}_i &=& v_i(v_i-\alpha) - w_i + gs(er-v_i) \\
\dot{w}_i &= &a(bv_i-w)\\
s &=& \frac{\bar{w}}{\eta} = \frac{1}{\eta} \left(\frac{1}{ N} \sum_{i=1}^N w_i    \right)\\
v_i(t_{spike}^-) &=& v_{peak} \rightarrow 
\begin{array}{rcl}
v_i(t_{spike}^+) &=& v_{reset} \\
w_i(t_{spike}^+) &=& w_i(t_{spike}^-) + \hat{w}, 
\end{array}
\end{eqnarray*}
for $i=1,2,\ldots N$. Here the dynamics of $s$ are replaced entirely by its steady state, large network solution: $\tau_s s_{jump} \langle R \rangle \approx w \frac{s_{jump}{\tau_s}}{w_{jump}\tau_s} = w /\eta$, with $w$ replaced by the finite mean $\bar{w}$.   We will refer to this network of neurons as the slow network.  

\begin{figure}
       \centering

         \subfigure[$g = g^*- 1$]{\includegraphics[scale=0.51]{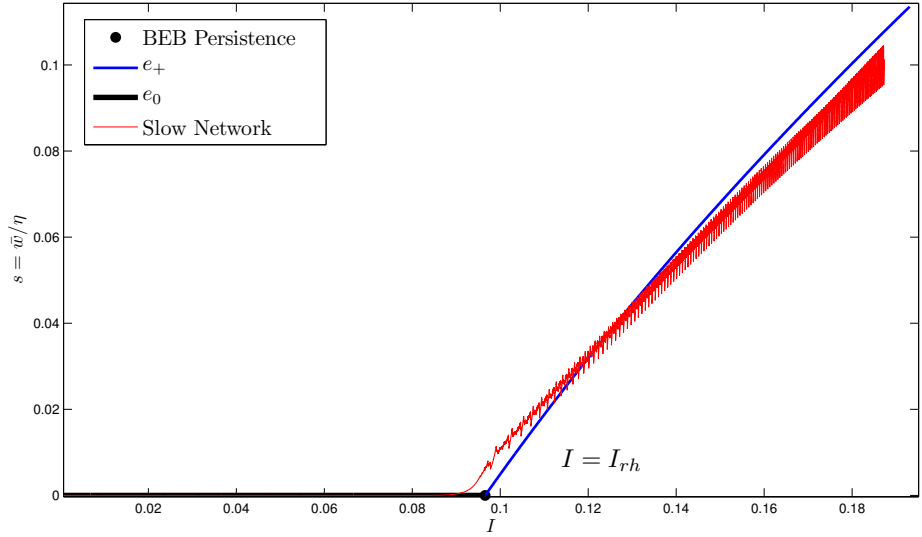}
\label{fig11a}}
\\
         \subfigure[$g = g^*+1$]{\includegraphics[scale=0.51]{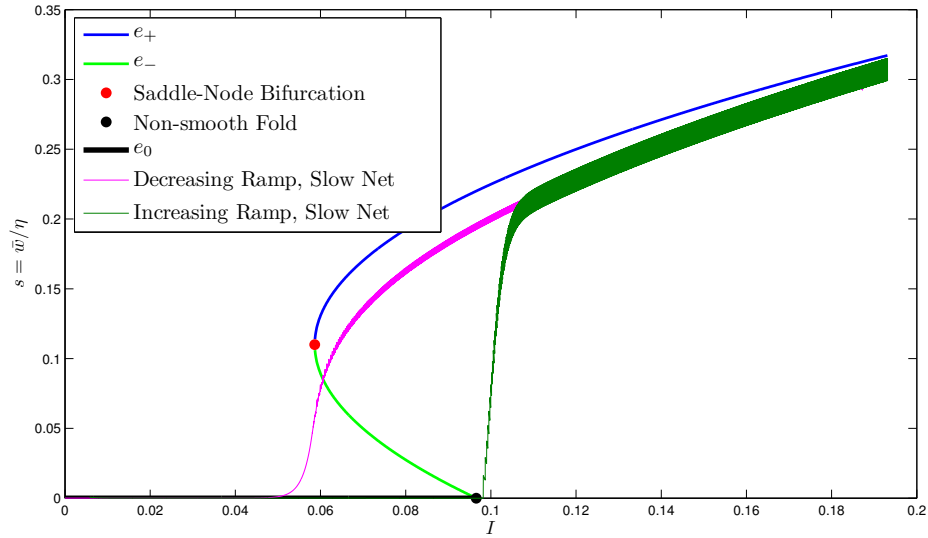}
\label{fig11b}}

\caption{The slow network consisting of 100 neurons is simulated with a slow current ramp.  The current is either descending (red) or ascending (green).  When $g<g^*-1$, no bistability is present and the steady state solution for the network collides with the non-firing solution, as predicted by the mean-field analysis.  When $g>g^*$, the descending current results in firing for $I<I_{rh}$, until the steady state falls off sharply near $I=I_{SN}$, as predicted by the mean-field analysis.  The ascending current only results in firing when $I=I_{rh}$ is reached.  
As there is bistability between these two stable states for the network for $g>g^*$, we should expect an unstable steady state separating these two, which is what the mean-field system predicts.} \label{fig11}
\end{figure}

The resulting mean-field system for the slow network is simply a one-dimensional non-smooth ODE, given by:
\begin{eqnarray*}
\dot{w} &=& -\frac{w }{\tau_w} + w_{jump}\langle R \rangle \\
\langle R \rangle &=& \begin{cases} \sqrt{ I-I^*(w /\eta,w)}& I \geq I^*(w /\eta,w)\\ 0 &  I<I^*(w /\eta,w)  \end{cases}
\end{eqnarray*}
The mean-field system for the slow network has the same steady states as the mean-field system for the full network, the two solutions $w_\pm$, in addition to the non firing solution, $w_0 = 0$.  However, being a one dimensional system, no Hopf bifurcations (and thus oscillations) are present in the mean-field system for the slow network.  Additionally, one can show that $w_+$ is always stable, and $w_-$ is always unstable as the eigenvalue for the steady state is given by 
\begin{equation}
\lambda(w_\pm) = -\frac{\lambda_w^2}{\tau_w}A_2(g) \left(1 - \frac{M(g)(g-g^*)}{M(g)(g-g^*)\pm \sqrt{M(g)^2(g-g^*)^2+\tilde I}}\right).
\end{equation}
The functions $A_2(g)$, and $M(g)$ are identical to those for the mean-field system for the full network, thus we should expect that the slow network has a stable steady state, and undergoes a saddle node bifurcation as for the full network.   As there can be no oscillations, we expect that $w_+$ exists and is stable for $g<g^*$, and $I>I_{rh}$ and $g>g^*$, and $I>I_{SN}$.  As the non-firing solution is also stable for $I<I_{rh}$, then we should expect bistability for $I_{SN}<I<I_{rh}$.  Indeed, if we simulate the slow network with a slowly varying current that either decreases from current values greater then $I_{rh}$ or increases from current values less then $I_{rh}$, we get bistability for $g>g^*$ and none for $g<g^*$. This is shown in Figure \ref{fig11} 

\begin{figure}
        \centering
         \subfigure[$g = g^*- 1$]{\includegraphics[scale=0.1]{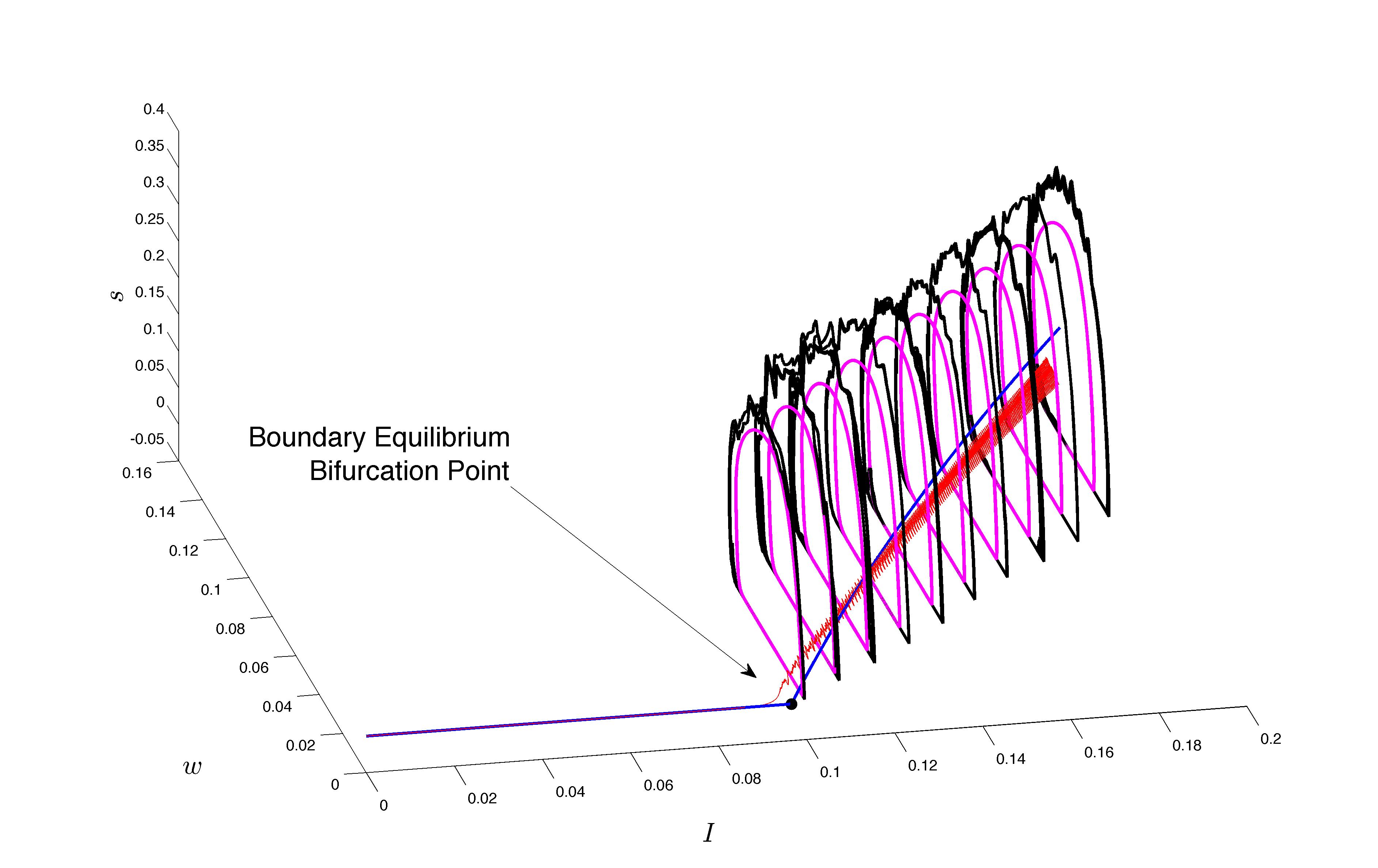}
\label{fig12a}}
\\
         \subfigure[$g = g^*+1$]{\includegraphics[scale=0.1]{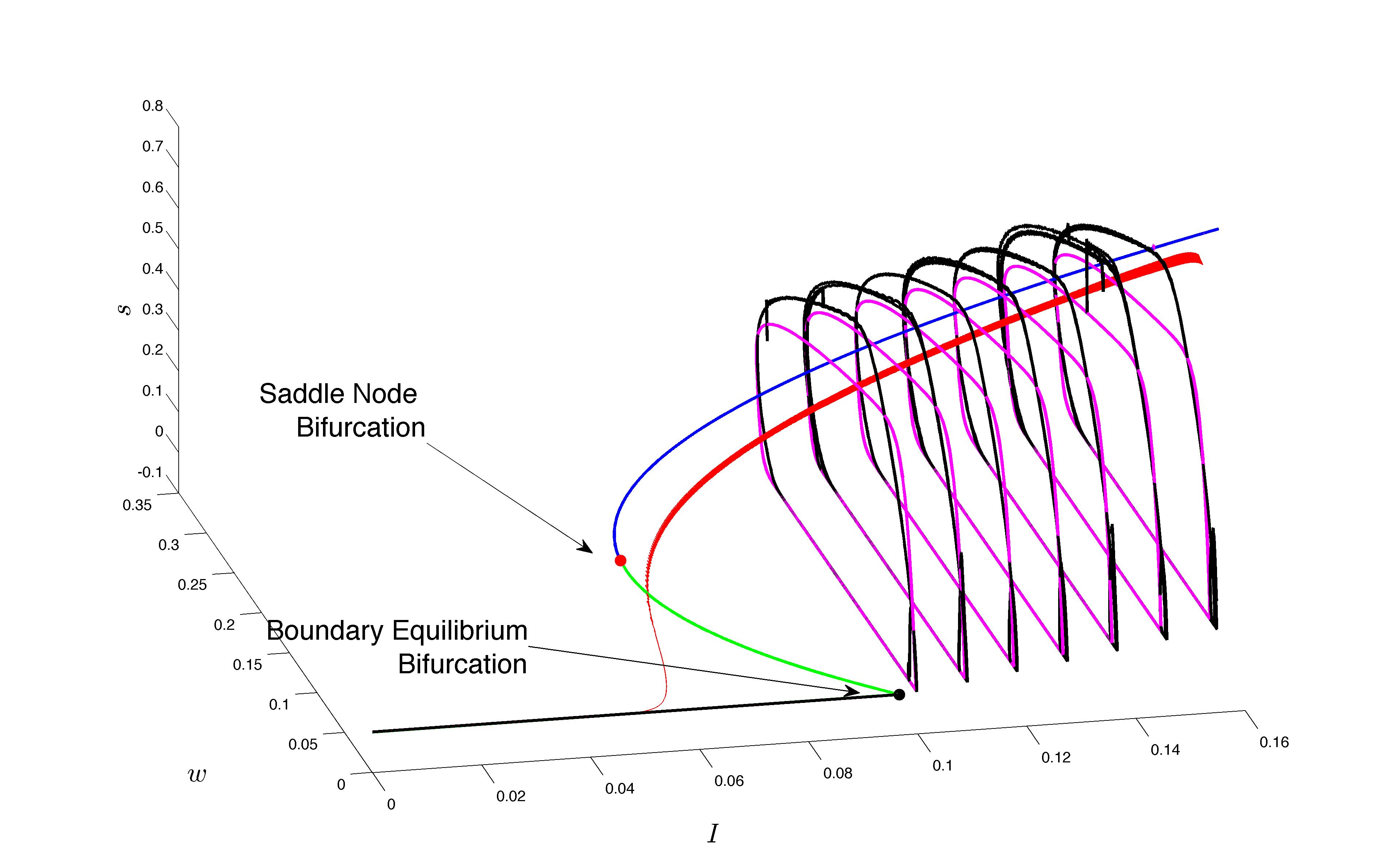}
\label{fig12b}}
\caption{Using simulations of the slow network to converge to the non-bursting steady state, and the full network to converge to the stable bursting limit cycle, we can piece together a pseudo-bifurcation diagram for the full network of neurons that very closely mirrors the bifurcation diagram predicted from the non-smooth mean-field equations.  Indeed, it appears that the transitions that occur at $I=I_{rh}$ are well explained as non-smooth boundary equilibrium bifurcations for the full network.   This would seem to indicate the existence of the co-dimension 2 non-smooth saddle-node BEB point for the actual network as well.  Note that the limit cycles have been smoothed out somewhat for clarity in the image, removing some of the high frequency oscillations due to synchrony in the peaks.} \label{fig12}
\end{figure}

Using the simulations of the slow network and the full network, we can piece together a pseudo-bifurcation diagram for the full network.  This is shown in figure \ref{fig12}.  The boundary equilibrium bifurcations that occur near the vicinity of $g^*$ are both present for the actual network. Given the similarities between the bifurcation diagram for the actual network, and that predicted by the non-smooth mean-field equations, it would appear that in order to understand the bifurcations that occur for these networks, one has to consider non-smooth bifurcation theory.   

One might wonder as to whether or not the non-smooth nature of the mean-field system is a direct result of the non-smooth nature of the neurons, given the fact that they have discrete resets and jumps.  However, this is not the case.  In particular, the firing rate of all type I neurons in the vicinity of the saddle-node on an invariant circle bifurcation is always proportional to $\sqrt{I-I_{rh}}$.  Assuming that the dynamics of the neurons voltage is much faster then the dynamics of the all the other intrinsic and synaptic currents, one could obtain mean-field equations very similar to that obtained explicitly here.  For example, this is done in the finite network case in the work of \cite{bardpaper}.   Thus, one has to consider non-smooth bifurcations and bifurcation analysis when working with mean-field systems for type I neurons.  The same is true of type 2 neurons, however the firing rate for these neurons changes discontinuously at $I_{rh}$, and thus it is likely that the mean-field systems for type 2 neurons would be completely non-smooth, as opposed to piecewise smooth continuous.

\section{Discussion}

Through our analysis of the mean field systems for large networks of coupled neurons, a number of new non-smooth bifurcations have been discovered that have been previously been unknown in the literature.  These include two co-dimension 1 branches of boundary equilibrium bifurcations that have homoclinic limit cycles at the bifurcation point, and can be thought of as generating/destroying non-smooth limit cycles.  Additionally, a pair of co-dimension 2 bifurcations have also been discovered that result from the collision of classical smooth branches of bifurcations with non-smooth bifurcations.  These occur when either a Hopf equilibrium point, or a saddle-node equilibrium point collide with a switching manifold.  We have determined locally in a neighbourhood of these bifurcation points the resulting behavior of the system through analytical and numerical results.  Additionally, a global-codimension 2 bifurcation and the collision of a Bogdanov-Takens equilibrium point with a switching manifold was also discovered by analyzing these systems, however we leave their analysis for later work.   

Given the analysis we have performed, then there are several predictions we can make, using the bifurcation curves for all the necessary branches.  For example, it appears that the time scales $\tau_w$ and $\tau_s$ are crucial for determining the presence of bursting.  If $\tau_w<\tau_s$, then no bursting can occur, while if $\tau_w>\tau_s$, there is a bell shaped region of bursting for $I>I_{rh}$, and $g>\bar{g}$.   Thus, if the adaptation time scale is smaller then the time scale of the synapses, adapting, recurrently coupled networks would not burst.   This can be the case for example for weakly adapting neurons coupled together with NMDA synapses, which have a long time scale.

Unfortunately however, while much of these bifurcations can be at least derived, one cannot easily determine whether or not even the smooth bifurcations are generic in any sense.  This is due to the fact that these particular systems cannot be diagonalized very easily, due to the presence of the $\sqrt{I-I^* (s,w)}$ term, which has an unbounded derivative as $I\rightarrow I^*(s,w)$.  Thus, center manifold theory cannot be directly applied, and many of the genericity conditions cannot be verified.   

In addition to these problems with regards to smooth bifurcations, one cannot easily apply the existing non-smooth theory to these systems.  The systems in the form 
\begin{eqnarray}
\dot{s} &=& f(s,u)= -\frac{s}{\tau_s} + \hat{s}\langle R_i(t) \rangle  \label{s2}\\
\dot{w} &=& g(s,u)=-\frac{w}{\tau_w} + {\hat{w}}\langle R_i(t) \rangle \label{w2}\\
\langle R_i(t) \rangle &=& \begin{cases} \sqrt{I-I^*(s,w )} & \quad  I\geq I^*(s,w )\\  0 & \quad  I<I^*(s,w )\end{cases} 
\end{eqnarray}
are clearly piecewise smooth continuous, however unlike the vast majority of PWSC systems discussed in the literature, they fail to satisfy one critical constraint that these other systems have.  In normal piecewise smooth continuous systems, given by 
\begin{eqnarray*}
\dot{x} = \begin{cases} f_1(x) & \quad \text{if} \quad H(x)\geq 0 \\ f_2(x) &\quad \text{if} \quad H(x)<0 \end{cases}
\end{eqnarray*}
where $f_1(x) = f_2(x)$ on $H(x)=0$, it is assumed that both $f_1(x)$ and $f_2(x)$ exists everywhere, and are smooth.  In our system, $\sqrt{I-I^* (s,w)}$ only exists when $I\geq I^*(s,w)$ and its first derivative only exists when $I>I^*(s,w)$.  This renders much of the analysis on PWSC systems inapplicable.  Indeed, this system cannot even be regularized in the normal way, via a Teixeira type regularization scheme \cite{Tex} due to the fact that $\sqrt{I-I^*(s,w)}$ is undefined when $I<I^*(s,w)$.  

However, there are alternate ways to apply both center manifold theory, and simultaneously regularize this system.   In particular, consider the three-dimensional system given by 
\begin{eqnarray*}
\dot{s} &=& -\frac{s}{\tau_s} + \hat{s}R \\
\dot{w} &=& -\frac{w}{\tau_w} + \hat{w}R \\
\epsilon \dot{R} &=& R(R^2 - (I-I^*(s,w)))
\end{eqnarray*}
where $\epsilon$ is a small constant.  In this singularly perturbed system, one can show that as $\epsilon \rightarrow 0$, one recovers the piecewise smooth continuous system as when $\epsilon$ is small, we can regard $s$ and $w$ as fixed, and thus $R$ rapidly converges to the steady states $0$, or $\sqrt{I-I^*(s,w)}$, depending on the sign of ${I-I^*(s,w)}$.  

Using this type of embedded system, which is entirely polynomial for the Izhikevich mean field system, one has regularized the non-smooth system in a sense by embedding it as the fast system in a singular perturbation problem.   Thus, to actually determine the genericity properties of the bifurcations displayed above, in addition to how the non-smooth bifurcations discovered are related to the general smooth bifurcation theory, one can analyze either directly or numerically the embedded system for finite $\epsilon$.  

Preliminarily, we have found that the Hopf BEB bifurcation seems to be a Bautin point under the regularization, and the saddle-node BEB bifurcation point seems to be a Bogdanov-Takens point under the regularization.  These are both co-dimension 2 smooth bifurcation, and they also explain the emergence of the non-smooth saddle-node of periodics in the Hopf BEB, which has a smooth saddle-node of periodics.  Additionally, the regularized Bogdanov-Takens has (generically) a branch of homoclinic bifurcations,  which also exists in a non-smooth form for the saddle-node BEB.   However, as the embedded system, and its justification as a regularization are outside of the scope of this paper, we leave it for future work.

In addition to the embedded regularization, this system is also unusual in the sense that there is a natural regularization for the mean-field system.  Suppose we consider the voltage equations to be perturbed by white noise:
\begin{equation}
\dot{v}_i = v_i(v_i-\alpha) - w_i + gs(er-v_i) + I + \eta_i 
\end{equation}
where $\langle \eta_i(t)\rangle = 0$ and $\langle \eta_i(t)\eta_i(t')\rangle = \sigma^2\delta(t-t')$.  In which case one can rigorously derive a mean-field system for this network of equations which is identical to the original mean-field system given in equations (\ref{MFa})-(\ref{MFc}) only the firing rate is now given by:
\begin{eqnarray}
\langle R \rangle = \left[\int_{v_{reset}}^{v_{peak}}\int_{v'}^{v_{peak}} \exp\left(-\frac{2}{\sigma^2}(M(v',w,s)-M(v,w,s)\right)\,dv'dv \right]^{-1}
\end{eqnarray}
where $M(v,w,s)$ is an anti-derivative (in $v$) of $F(v) - w + gs(er-v) + I$.    As we shall do in forthcoming work \cite{us5} one can rigorously show that this expression for $\langle R\rangle$ is smooth with respect to $s$ and $w$, and always defined, and converges to (\ref{MFc}) as $\sigma \rightarrow 0$.  Thus, the mean-field system for a network with noise parameterized by $\sigma^2$, the variance in the noise, converges to the system (\ref{MFa})-(\ref{MFc}) as the variance of the noise becomes negligible.  But, since the mean-field system with noise is smooth, it can be thought of as a natural regularization for the non-smooth mean-field system.  We remark that this is unusual in the field of non-smooth theory as generally a regularization is chosen or suggested, and is typically of the Teixeira form \cite{Tex}.  We leave the bifurcation analysis of the system with noise for future work.  

Finally, one may ask if the non-smooth bifurcations we analyze here appear in other non-smooth systems or are generic in any way.  To the best of our knowledge, the co-dimension 2 bifurcations are novel in the literature, according to a recent review \cite{Jeffrey}.  However, a Hopf-bifurcation occuring on a discontinuity boundary (a co-dimension 2 non-smooth bifurcation) does occur in the example (in section 6) in \cite{Kuz} (see Figure 29.).  However, the system examined in \cite{Kuz} is a Fillipov system, and thus has a higher order of discontinuity.  

As to whether these bifurcations occur in a more generic system, we intend to explore this further with a more generic piecewise smooth continuous system that does not have undefined derivatives on the switching manifolds.  In particular, we note that with a quadratic PWSC system, one can show that as for our saddle-node BEB bifurcation, at the intersection between generic branches of persistence and non-smooth fold BEB bifurcations, the equilibrium of $f_2(x)$ must have a zero eigenvalue at this co-dimension 2 point.    We intend to explore all the possible cases and their possible relationships to the co-dimension 2 non-smooth bifurcations we outline in this paper for future work.


%
%
%


\begin{table}
\centering
\begin{tabular}{|c|c|c|c|}
\hline
Parameter & Izhikevich Network (from \cite{us} )& AdEx Network  (from \cite{AdEx}) & QIF Network (from \cite{Touboul2008} ) \\
\hline
$g$ & 0-4 & 0-1000 & 0-40\\
\hline
$I$ & 0-0.4 & 0 -60 & 0-40\\
\hline
$\tau_s$ & 2.6 & 0.08& 2\\
\hline
$\tau_w$ & 130 & 3.63&50\\
\hline
$s_{jump}$ & 0.8 &0.5&1\\
\hline
$w_{jump}$ & 0.0189& 21.92&0.36\\
\hline
$e_r$ & 1 & 2.5 & 2\\
\hline
$\alpha$ & 0.62& N/A & N/A \\
\hline
$v_{peak}$ & 1.46 & 65& 10\\
\hline
$v_{reset}$ & 0.15 & -1.25 & 0\\
\hline
\end{tabular}
\caption{Parameters for various network types and the mean-field systems.  Note that parameter above are dimensionless, where as in some of the cited sources they are in dimensional form only. } \label{table1}
\end{table}


\begin{thebibliography}{10}




%


\bibitem{AdEx}
R.~Brette and W.~Gerstner.
\newblock Adaptive exponential integrate-and-fire model as an effective
  description of neuronal activity.
\newblock {\em Journal of Neurophysiology}, 94(5):3637--3642, 2005.
\bibitem{Tex}
A.~Buzzi, P.R. da Silva, and M.A.~Teixeira  
\newblock A singular approach to discontinuous vector fields on the plane 
\newblock {\em Journal of Differential Equations}, 231:633--655, 2006

\bibitem{Jeffrey}
A.~Colombo, M.~di Bernardo, and M.R.~Jeffrey 
\newblock Bifurcations of piecewise smooth flows:  perspectives, methodologies and open problems 
\newblock {\em Physica D}, 241(22):1845--1860, 2011


\bibitem{nonsmooth}
M.~di~Bernardo, C.J. Budd, A.R. Champneys, P.~Kowalczyk, A.B. Nordmark, G.O.
  Tost, and P.T. Piiroinen.
\newblock Bifurcations in nonsmooth dynamical systems.
\newblock {\em SIAM Review}, 50(4):629--701, 2008.


\bibitem{MATCONT}
A.~Dhooge, W.~Govaerts, and Yu.~A. Kuznetsov.
\newblock {MATCONT}: A {MATLAB} package for numerical bifurcation analysis of {ODEs}.
\newblock {\em ACM Transactions on Mathematical Software}, 29(2):141--164, 2003.



\bibitem{us}
M.~{Dur-e-Ahmad}, W.~Nicola, S.A. Campbell, and F.K.~Skinner.
\newblock Network bursting using experimentally constrained single compartment
  {CA3} hippocampal neuron models with adaptation.
\newblock {\em Journal of Computational Neuroscience}, 33(1): 21--40, 2012.



\bibitem{bardbook}
G.B.~Ermentrout and D.H.~Terman 
\newblock Mathematical Foundations of Neuroscience
\newblock {\em Springer, New York}, 2010 

\bibitem{bardpaper}
G.B.~Ermentrout 
\newblock Reduction of conductance-based models with slow synapses to neural nets.
\newblock {\em Neural Computation}, 6(4):679--695 (1994) 




%





\bibitem{Katie}
K.A.~Ferguson, C.Y.L.~Huh, B.~Amilhon, S.~Williams, and F.K.~Skinner
\newblock Simple, biologically-constrained CA1 pyramidal cell models using an intact, whole hippocampus context
\newblock {\em F1000Research}, 2014, 3:104 (doi: 10.12688/f1000research.3894.1)






%
%
%

\bibitem{Izhikevich}
E.M. Izhikevich.
\newblock Simple model of spiking neurons.
\newblock {\em Neural Networks, IEEE Transactions on}, 14(6):1569 -- 1572, 2003.



%



\bibitem{Kuz}
Y.A..~Kuznetsov, S. Rinaldi, and A.~Gragnani  
\newblock One-parameter bifurcations in planar fillipov systems
\newblock {\em International Joural of Bifurcation and Chaos}, 13(8):2157--2188, 2003

\bibitem{Ad2}
R.~Naud, N.~Marcille, C.~Clopath, and W.~Gerstner.
\newblock Firing patterns in the adaptive exponential integrate-and-fire model.
\newblock {\em Biological Cybernetics}, 99:335--347, 2008.


\bibitem{us5}
W.~Nicola, C.~Ly and S.A.~Campbell.
\newblock One-dimensional population density approaches to recurrently coupled networks of neurons with noise 
\newblock {\em In preparation} 2014
%


\bibitem{us2}
W.~Nicola, and S.A.~Campbell.
\newblock Bifurcations of large networks of two-dimensional integrate and fire neurons 
\newblock {\em Journal of Computational Neuroscience}, 35: 87--108, 2013. 
%
%
\bibitem{NT}
D.Q. Nykamp. and D.~Tranchina.
\newblock A population density approach that facilitates large-scale modeling
  of neural networks: Analysis and an application to orientation tuning.
\newblock {\em Journal of Computational Neuroscience}, 8:19--50, 2000.
%
\bibitem{OKS}
A.~Omurtag, B.W. Knight, and L.~Sirovich.
\newblock On the simulation of large populations of neurons.
\newblock {\em Journal of Computational Neuroscience}, 8:51--63, 2000.

%
%
%
%

\bibitem{Touboul2008}
J.~Touboul.
\newblock Bifurcation analysis of a general class of nonlinear
  integrate-and-fire neurons.
\newblock {\em SIAM Journal on Applied Mathematics}, 68(4):1045--1079, 2008.









%
%
%
\bibitem{VreeswijkandHansel}
C.~van Vreeswijk and D.~Hansel.
\newblock Patterns of synchrony in neural networks with spike adaptation.
\newblock {\em Neural Computation}, 13(5):959--992, 2001.
%
%


\end{thebibliography}
\end{document}